\documentclass[12pt]{article}

\usepackage{graphicx}  
\usepackage{amsmath}  
\usepackage{mathtools}
\usepackage{amssymb,amsthm}
\usepackage{color}
\usepackage{tikz}
\usepackage{pgfplots}
\usepackage[pdfpagemode=UseOutlines]{hyperref}
\hypersetup{
    colorlinks=false,
    pdfborder={0 0 0},
    pdfauthor={T. Haas, B. de Rijk, G. Schneider},
   pdftitle={Modulation equations near the  Eckhaus boundary: the KdV equation}
}
\pgfplotsset{compat=1.11}

\newcommand{\red}{\color{red}}

\newcommand{\R}{{\mathbb R}}\newcommand{\N}{{\mathbb N}}
\newcommand{\Res}{\mathrm{Res}}
\newcommand{\C}{{\mathbb C}}
\newcommand{\e}{e}
\renewcommand{\i}{\mathrm{i}}

\newcommand{\fourier}{\mathcal{F}}

\def\per{\mathrm{per}}
\def\F{\mathcal{F}}

\def\A{\mathcal{A}}
\def\cR{\mathcal{R}}

\def\h{\mathcal{H}}

\def\E{\mathcal{E}}
\def\Non{\mathcal{N}}
\def\bd{\mathrm{bd}}

\def\El{\mathcal{L}}
\def\Re{\mathrm{Re}}
\def\Im{\mathrm{Im}}
\def\ord{\mathcal{O}}

\def\Y{\mathcal{Y}}
\def\cW{\mathcal{W}}
\def\cZ{\mathcal{Z}}

\let\tht\theta
\let\epsilon\varepsilon
\let\theta\vartheta
\let\hat\widehat

\newtheorem{theorem}{Theorem}[section]\newtheorem{lemma}[theorem]{Lemma}

\newtheorem{corollary}[theorem]{Corollary}

\newtheorem{remark}[theorem]{Remark}

\title{Modulation equations near \\ the  Eckhaus boundary \\
-- The KdV equation --}
\author{Tobias Haas\footnotemark[2] \footnotemark[4]
\and Bj\"orn de Rijk\footnotemark[2] \footnotemark[3] 
\and Guido Schneider\footnotemark[2] \footnotemark[5]}
\date{\vspace{-2em}}

\begin{document}

\maketitle

\renewcommand{\thefootnote}{\fnsymbol{footnote}}
\footnotetext[2]{Institut f\"ur Analysis, Dynamik und Modellierung, Universit\"at Stuttgart, Pfaffenwaldring 57, 70569 Stuttgart, Germany}
\footnotetext[4]{tobias.haas@mathematik.uni-stuttgart.de}
\footnotetext[3]{bjoern.derijk@mathematik.uni-stuttgart.de}
\footnotetext[5]{guido.schneider@mathematik.uni-stuttgart.de}

\begin{abstract}
We are interested in the description of small modulations in time and space of wave-train solutions to the complex Ginzburg-Landau equation
\begin{align*}
\partial_T \Psi = (1+ i \alpha) \partial_X^2 \Psi + \Psi - (1+i \beta ) \Psi |\Psi|^2
\end{align*}
near the Eckhaus boundary, that is, when the wave train is near the threshold of its first instability. Depending  on the parameters $ \alpha $, $ \beta $ a number of modulation equations can be derived, such as the KdV equation, the Cahn-Hilliard equation, and a family of Ginzburg-Landau based amplitude equations. Here we establish error estimates showing that the KdV approximation makes correct predictions in a certain parameter regime. Our proof is based on energy estimates and exploits the conservation law structure of the critical mode. In order to improve linear damping we work in spaces of analytic functions.

\textbf{Keywords.} Modulation equation, validity, wave trains, long wave approximation, Eckhaus boundary
\end{abstract}



\newpage 

\section{Introduction}

The complex Ginzburg-Landau (GL) equation
\begin{equation} \label{eq1}
\partial_T \Psi = (1+ i \alpha) \partial_X^2 \Psi + \Psi - (1+i \beta ) \Psi |\Psi|^2, \quad X \in \mathbb{R}, T \geq 0, \Psi(X,T) \in \mathbb{C},
\end{equation}
with $\alpha,\beta \in \R$ can be derived via a multiple scaling analysis as a universal amplitude equation for the description of pattern forming systems, such as reaction-diffusion systems or the Couette-Taylor problem, near the threshold of the first instability of the trivial ground state, cf.~\cite{NW69}. See~\cite[Chapter 10]{SUbook} for a recent survey.

The GL equation possesses a family of wave-train solutions
\begin{equation} \label{equi1}
\Psi_{\per}(X,T) = \Psi_0 e^{i(\zeta X + \Omega_0 T)}, 	
\end{equation}
which are periodic in time and space, with $\Psi_0 , \Omega_0 \in \R$ and $\zeta \in (-1,1)$ satisfying
\begin{align}\Psi_0^2 + \zeta^2 = 1,  \qquad \Omega_0 + \alpha \zeta^2 + \beta \Psi_0^2 =0. \label{existcond}\end{align}
As~\eqref{eq1} is invariant under the mapping $\Psi \mapsto - \Psi$, we assume without loss of generality $\Psi_0 > 0$ throughout the paper. Moreover, spatially homogeneous wave trains lie outside the parameter regime -- see~\S\ref{sec2} -- considered in this paper and so we assume $\zeta \in (-1,1) \setminus \{0\}$ in the following.

The stability of these solutions was first discussed in~\cite{Eck65}. Due to translational invariance of the family of wave trains in time and space as  solutions to~\eqref{eq1}, the spectrum of the linearization of~\eqref{eq1} about $\Psi_{\per}$ touches the origin. Therefore, in the most stable scenario, the spectrum is bounded away from the imaginary axis in the left-half plane except for a  tangency at the origin. In such a case, we call the wave train {spectrally stable}.

It is well-known~\cite{Eck65}, that in case $\alpha = \beta = 0$ wave trains are spectrally stable if and only if $\zeta^2 \leq 1/3$. In fact, for all fixed $\alpha, \beta \in \R$ with $1 + \alpha\beta > 0$, there is a critical wave number $\zeta_{\bd} = \zeta_{\bd}(\alpha,\beta) \in (0,1)$ -- the so-called {Eckhaus boundary} -- such that spectral stability holds if and only if $|\zeta| \leq  \zeta_{\bd}$, cf.~\cite{vH94}. We note that in case $|\zeta| <  \zeta_{\bd}$ spectral stability yields nonlinear stability of wave-train solutions to~\eqref{eq1} with respect to small spatially localized perturbations~\cite{BK92,CEE92,Kap94}. The same result at the Eckhaus boundary $|\zeta| = \zeta_{\bd}$ for $ \alpha =  \beta = 0 $ has been  established in~\cite{GSWZ17}.

\begin{figure}[b!]
 \centering
 \includegraphics[width=0.5\linewidth]{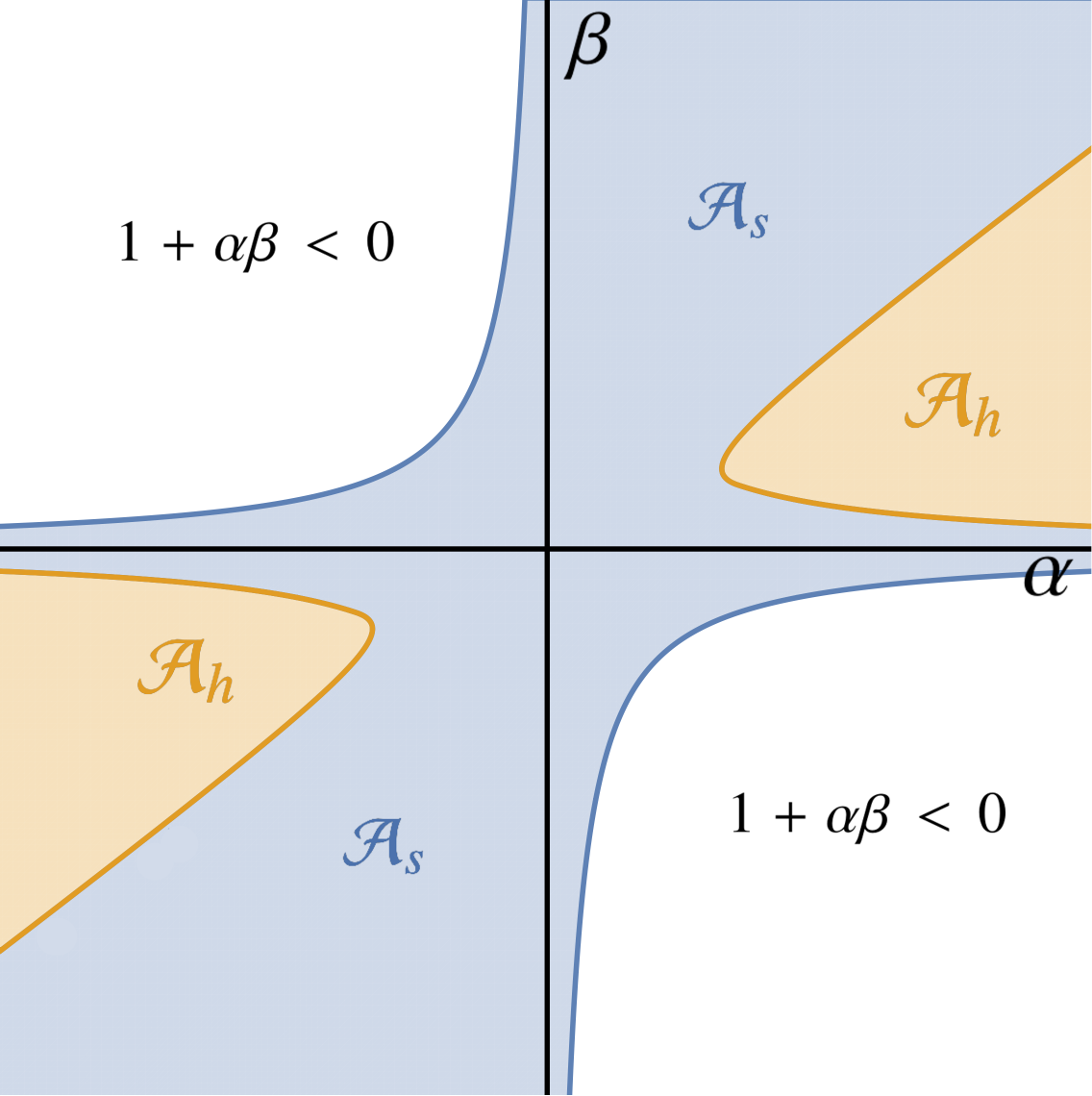}
 \caption{The regions $\mathcal{A}_s$ and $\mathcal{A}_h$ in the $(\alpha,\beta)$-plane.}
\label{vanHarten}
\end{figure}

It depends on the value of $(\alpha,\beta) \in \R^2$ whether the wave train $\Psi_\per$ destabilizes through a Hopf-Turing or sideband instability at the Eckhaus boundary $|\zeta| = \zeta_{\bd}$. More specifically, there are disjoint open regions $\A_s, \A_h \subset \R^2$ (which can be determined explicitly~\cite{vH94}) satisfying $\overline{\A_s \cup \A_h} = \{\alpha,\beta \in \R^2 : 1+\alpha\beta \geq 0\}$ such that a Hopf-Turing instability occurs at $|\zeta| = \zeta_{\bd}$ if $(\alpha,\beta) \in {\A_h}$ and a sideband instability occurs at $|\zeta| = \zeta_{\bd}$ if $(\alpha,\beta) \in{\A_s}$ -- see also Figure~\ref{vanHarten}.

In~\cite{vH94} various amplitude equations for the description of slow modulations in time and space of the wave-train solutions to~\eqref{eq1} have been formally derived for $|\zeta|$ close to the Eckhaus boundary $\zeta_{\bd}$. One obtains a  system consisting of a Ginzburg-Landau equation
coupled to a nonlinear diffusion equation
 if $(\alpha,\beta) \in \A_h$, a Korteweg-de Vries (KdV) equation if $(\alpha,\beta) \in \A_s$ with $\alpha \neq \beta$ and a Cahn-Hilliard equation if $(\alpha,\beta) \in \A_s$ with $\alpha = \beta$.
At the boundaries of $ \A_h $ and $ \A_s $ more complicated amplitude equations occur.

In the last decades it has been observed that the Eckhaus boundary plays an important role in the
creation of patterns, especially in the wave number selection of the pattern, cf.~\cite{ahlers}.
Thus, for an analytic understanding of these pattern forming processes,
it is important to know which of the aforementioned amplitude equations occurring at the Eckhaus boundary
are valid and which are not.
\newpage
It is the purpose of this paper to establish error estimates showing that the KdV equation
\begin{align} \partial_{\tau} A  = (\beta - \alpha)\left(\frac{1+\alpha\beta}{1+\beta^2} \ \partial_{\xi}^3 A  +  \partial_{\xi}(A^2)\right), \quad \xi \in \R, \tau \geq 0, A(\xi,\tau) \in \mathbb{R}, \label{KdV}\end{align}
makes correct predictions about the dynamics of slow modulations of the wave-train solution $\Psi_{\per}$ to~\eqref{eq1}. Thus, we assume $(\alpha,\beta) \in \A_s$ with $\alpha \neq \beta$ and $ |\zeta| = \zeta_{\bd} + \ord(\epsilon^2)$ -- where $0 < \epsilon \ll 1$ is a small perturbation parameter -- such that
for $  |\zeta| > \zeta_{\bd} $
a small unstable branch of spectrum has entered the right-half plane after a sideband instability -- see Figure~\ref{sidebandfig}. In this regime the underlying structure driving the slow modulations of the wave train in the dissipative Ginzburg-Landau equation is conservative and governed by the KdV equation. We exploit this conservative structure to obtain nontrivial error estimates on an $\mathcal{O}(1/\varepsilon^3)$-time scale for solutions of order $\mathcal{O}(\varepsilon^2)$. Consequently, we do not only employ dissipative methods, but rely instead on analytical smoothing to improve linear damping in our error estimates. More precisely, we split the equations for the error into a linearly  exponentially damped part and into the rest. For estimating the rest we exploit improved linear damping due to analytical smoothing and the fact that the associated nonlinear term is $\mathcal{O}(k)$ for $ k \rightarrow 0 $.

Other rigorous  approximation results exist away from the Eckhaus boundary. In case $\alpha= \beta = 0$
and  $|\zeta| < \zeta_{bd}$
the so-called phase-diffusion equation has been justified~\cite{MS04a}. For $(\alpha,\beta) \neq (0,0)$ the validity of a conservation law~\cite{MS04b} and, again  for $|\zeta| < \zeta_{\bd}$, of the Burgers equation~\cite{DSSS09} has been established. On the other hand, {at} the Eckhaus boundary $|\zeta| = \zeta_{\bd}$ in the regime $\alpha = \beta = 0$, it is shown in~\cite{DS09JNS} that a waiting time phenomenon occurs.

A theorem about  the KdV approximation at the Eckhaus boundary
has already been stated in~\cite[\S 7.5]{DSSS09}. However, no detailed proof was given.
It was outlined how the proof for the validity of Burgers approximation would transfer to the KdV approximation at the Eckhaus boundary.
When
preparing the manuscript~\cite{dRHS18} we recognized that a complete validity proof
is much more involved   and goes far beyond the sketch  given in~\cite[\S 7.5]{DSSS09}.
Therefore, we decided to give a complete proof
for the validity of the KdV approximation at the
Eckhaus boundary of which this paper is the outcome.

The plan of the paper is as follows. In~\S\ref{sec1b} we derive
equations for modulations of the wave train in a suitably chosen coordinate system. In~\S\ref{sec2}
we recall the calculations from~\cite{vH94} to determine the parameter regime which
leads to the region $\A_s$ in Figure~\ref{vanHarten} and to the derivation of the
KdV equation in~\S\ref{sec3}.
~\S\ref{setup} contains the functional analytic set-up.
In~\S\ref{mainresults} our approximation result is formulated.
The proof of this result is given in~\S\ref{secproof}.
\S\ref{secdiscussion} consists of concluding remarks in which we discuss the approximation result
in the original coordinate system, what happens if we work in Sobolev spaces, and
whether other formally derived amplitude equations at the Eckhaus boundary are valid or not.
Technical results are provided in an appendix.
\medskip

{\bf Notation.} Constants which can be chosen independently of the small
perturbation parameter $ 0 < \varepsilon \ll 1 $ are denoted by the same symbol $ C $.
\medskip

{\bf Acknowledgement.} We gratefully acknowledge financial support by the Deutsche
Forschungsgemeinschaft (DFG)
through CRC 1173.

\section{The equations for the modulation}
\label{sec1b}

We are interested in modulations of the wave-train solution~\eqref{equi1} to the GL equation~\eqref{eq1}. Thus, we consider the modulated solution
\begin{align} \Psi(X,T) = e^{r_0 + s(X,T) + i(\zeta X + \Omega_0 T + \phi(X,T))} = \Psi_{\per}(X,T) e^{s(X,T) + i\phi(X,T)}, \label{modsol}\end{align}
to~\eqref{eq1}, where $r_0 \in \R$ is defined through $\Psi_0 = e^{r_0}$, cf.~\eqref{equi1}. The conditions~\eqref{existcond} for the existence of wave-train solutions now read
\begin{align}
 e^{2 r_0} + \zeta^2 = 1 , \qquad  \Omega_0 + \alpha \zeta^2 + \beta  e^{2 r_0} =0. \label{existcond2}
\end{align}
In this section, we derive equations for the modulations $(s,\phi)$ in~\eqref{modsol}. Thus, we rewrite the GL equation~\eqref{eq1} in polar coordinates of the form $ \Psi = e^{r+ i \varphi} $, which yields
\begin{align}\begin{split}
\partial_T \varphi & = \partial_X^2 \varphi + \alpha \partial_X^2 r - \beta  e^{2r} + \alpha  (\partial_X r)^2- \alpha  (\partial_X \varphi)^2
+ 2 (\partial_X \varphi)(\partial_X r) , \\
\partial_T r & = \partial_X^2 r - \alpha \partial_X^2 \varphi + 1- e^{2r} + (\partial_X r)^2- (\partial_X \varphi)^2
-2 \alpha (\partial_X \varphi)(\partial_X r) ,
\end{split} \label{GLPC}\end{align}
where we used
\begin{align*}
\begin{split}
\partial_T \Psi &=  e^{r+ i \varphi}(\partial_T r+ i \partial_T \varphi), \\
\partial_X^2 \Psi &=  e^{r+ i \varphi}(\partial_X r+ i \partial_X \varphi)^2 + e^{r+ i \varphi}(\partial_X^2 r+ i \partial_X^2 \varphi) .
\end{split}
\end{align*}
The modulated solution~\eqref{modsol} in polar coordinates reads
\begin{align*} r(X,T) = r_0 + s(X,T), \qquad \varphi(X,T) = \phi(X,T) + \zeta X +  \Omega_0 T.\end{align*}
Substituting this solution into~\eqref{GLPC} and employing~\eqref{existcond2} yields equations for the modulations $U(X,T) = (\phi,s)(X,T)$, which are given by
\begin{align}
\partial_T U = \El_0 U + \Non_0 (U), \label{modeq1}
\end{align}
where $\El_0$ is the differential operator
\begin{align*}
\El_0 &= D_2 \partial_X^2 - 2\zeta D_1 \partial_X + e^{2r_0} D_0,
\end{align*}
with
\begin{align*}
D_2 = \left(\begin{array}{cc} 1 & \alpha \\ -\alpha & 1\end{array}\right), \qquad D_1 = \left(\begin{array}{cc} \alpha & -1 \\ 1 & \alpha\end{array}\right),\qquad D_0 = \left(\begin{array}{cc} 0 & -2\beta \\ 0 & -2 \end{array}\right),
\end{align*}
and the nonlinearity $\Non_0$ is given by
\begin{align*}
\Non_0(U) = \left(\begin{array}{c} \alpha (\partial_X s)^2 -\beta  e^{2r_0}h(s) - \alpha  (\partial_X \phi)^2 + 2 \partial_X \phi \partial_X s\\
  (\partial_X s)^2 -e^{2r_0}h(s) - (\partial_X \phi)^2 - 2 \alpha \partial_X \phi \partial_X s\end{array}\right),
\end{align*}
with $h(s) := e^{2s}-1-2s$. Following~\cite{vH94}, it is advantageous to switch to a co-moving frame and to rescale both time and space. Therefore, we introduce
\begin{align} \sigma = \frac{e^{2r_0}}{\zeta^2} = \zeta^{-2} - 1, \qquad t = \zeta^2 T,  \qquad x = \zeta (X - c \zeta T ),\label{newcoord}\end{align}
where we recall that   $|\zeta| \in (0,1)$. With respect to these coordinates, system~\eqref{modeq1} reads
\begin{align}
 \partial_t U = \El U + \Non(U),\label{modeq2}
\end{align}
where $\El$ denotes the differential operator
\begin{align*}
\El =  D_2 \partial_x^2 + (c I - 2 D_1) \partial_x + \sigma D_0,
\end{align*}
 and the nonlinearity $\Non$ is given by
\begin{align*}
 \Non(U) = \left(\begin{array}{c}  \alpha  (\partial_x s)^2-\sigma   \beta h(s) -  \alpha  (\partial_x \phi)^2 + 2 \partial_x \phi \partial_x s \\
 (\partial_x s)^2 -\sigma h(s) - (\partial_x \phi)^2 - 2  \alpha \partial_x \phi \partial_x s\end{array}\right).
\end{align*}
The modulation equation~\eqref{modeq2}  
only depend on $x$- and $t$-derivativ\-es of $\phi$ and not on $\phi$ itself. This yields translational invariance of the wave train $\Psi_{\per}$, i.e., any translate $\Psi_\per(X + X_0,T + T_0) = \Psi_\per(X,T) e^{i(\zeta X_0+\Omega_0T_0)}$, with $X_0,T_0 \in \R$,  of the wave-train solution in space and/or time is again a solution to~\eqref{eq1}. Consequently, the modulation equation~\eqref{modeq2}  admits constant solutions $(0,\phi_0)$ with $\phi_0 \in \R$. To account for translational invariance, it is beneficial to introduce the local {wave number} $\psi=\partial_x \phi $
of the modulated wave train. Then~\eqref{modeq2} transforms into
\begin{align} \partial_t V = L V + N(V),\label{modeq3}\end{align}
where $V(x,t) = (\psi,s)(x,t)$, where $L$ denotes the differential operator
\begin{align*}
L &= \left(\begin{array}{cc} \partial_x^2 + (c-2\alpha)\partial_x  & \alpha \partial_x^3 + 2\partial_x^2 - 2\sigma\beta\partial_x \\ -\alpha \partial_x - 2 & \partial_x^2 + (c-2\alpha)\partial_x - 2\sigma \end{array}\right),
\end{align*}
and the nonlinearity $N$ is given by
\begin{align}
 N(V) = \left(\begin{array}{c} \partial_x\left(\alpha (\partial_x s)^2 - \sigma \beta h(s) - \alpha  \psi^2 + 2 \psi \partial_x s \right) \\
 (\partial_x s)^2 - \sigma h(s)  - \psi^2 - 2  \alpha \psi \partial_x s \end{array}\right). \label{nonl1}\end{align}
Due to the introduction of the local wave number $\psi = \partial_x \phi$, we obtain a derivative $\partial_x$ in front of the first component of the nonlinearity $N$. Hence, we gain that the first component of the nonlinearity vanishes like $\mathcal{O}(k)$ in Fourier space at the wave number  $k = 0$. On the other hand, the introduction of $\psi$ yields a third derivative in the off-diagonal entry of the linearity $L$. This complicates establishing regularity and damping properties of $L$ in $L^2$-type spaces, which is required for our analysis in~\S\ref{secproof}. Therefore, we replace in~\S\ref{secproof} the derivative $\partial_x \phi$ in~\eqref{modeq2} by a {local} derivative $\theta(\phi)$ instead, which is defined through its action in Fourier space. It holds $\theta(\phi) = \F^{-1}[\hat{\theta} \F(\phi)]$ with $\hat{\theta}(k) = ik \min\{1,|k|^{-1}\}$, where $\F$ is the Fourier transform, cf. Figure~\ref{fig:pseudo-derivative}.

\section{Determining the parameter regime}
\label{sec2}

The aim of this paper is to prove that the KdV equation makes correct predictions about the behavior of modulations of marginally sideband-unstable wave-train solutions $\Psi_\per$ to~\eqref{eq1}. In this section, we establish a parameter regime that leads to the desired spectral configuration.

The linearization of~\eqref{eq1} about the wave-train solution $\Psi_\per$ in $(x,t)$-coordinates~\eqref{newcoord} corresponds to the linearity $\El$ in the modulation equation~\eqref{modeq2}. The spectrum of the operator $\El$ on $L^2(\R)$ is given by the eigenvalues of its Fourier symbol
\begin{align*}
\hat{\El}(k) =  \left( \begin{array}{cc} - k^2 + (c - 2 \alpha ) i k &
-  \alpha k^2 - 2 \beta \sigma + 2 ik
\\  \alpha k^2  - 2  i  k &  - k^2 - 2 \sigma + (c - 2 \alpha ) i k \end{array} \right).
\end{align*}
Equating $\det(\hat{\El}(k) - \lambda I) = 0$ gives
\begin{align*}\begin{split}
\tht (\tht - 2 \sigma )+\sigma^2 \gamma(k) ( \gamma(k) + 2 \beta  )
= 0,
\end{split}\end{align*}
where $ \tht = - \lambda - k^2 + (c - 2 \alpha ) i k $ and $\gamma(k) = (\alpha k^2 - 2 i k)/\sigma$. Solving this quadratic equation in $\tht$ yields two solutions
\begin{align*}\begin{split}
\tht_\pm & = \sigma \pm \sqrt{\sigma^2 -  \sigma^2 \gamma(k) ( \gamma(k) + 2 \beta  )} = \sigma \pm  \sigma \upsilon(k),
\end{split}\end{align*}
where the quantity $\upsilon(k) $ is defined as the principal square root
$$
v(k) = \sqrt{1-\gamma(k)^2-2 \beta \gamma(k)}.
$$
Consequently, the eigenvalues of $\hat{\El}(k)$ are
\begin{align}\begin{split}
\lambda_{\pm}(k) &\! = i(c-2\alpha)k - k^2 - \sigma \pm \sigma  \upsilon(k)\\
 &\!= i(c-2\alpha)k - k^2 - \sigma \pm \sqrt{\sigma^2 - \left(\alpha k^2 - 2 i k + 2 \beta\sigma\right)\left(\alpha k^2 - 2 i k\right)}.
\end{split}
\label{speccurv}
\end{align}
Thus, the spectrum of $\El$ is given by the union $\lambda_+[\R] \cup \lambda_-[\R]$. The curve $\lambda_-[\R]$ lies in the open left-half plane. On the other hand, $\lambda_+[\R]$ touches the origin at $k = 0$ implying $0$ is in the spectrum of $\El$, which is related to the translational invariance of the wave train $\Psi_\per$ in space and time as a solution to~\eqref{eq1}. We expand $\lambda_+$ about the origin as
\begin{align}
\lambda_{+}(k) & = c_1 i k - c_2 k^2 + c_3ik^3 - c_4k^4 + \mathcal{O}(|k|^5),\label{critdisp}
\end{align}
with
\begin{align*}\begin{split}
c_1 & = c-2(\alpha-\beta), \\
c_2 &= 1 + \alpha \beta - 2 (1+\beta^2) \sigma^{-1}, \\
c_3 &= 2(1+\beta^2)(\alpha\sigma - 2\beta)\sigma^{-2},\\
c_4 &= (1+\beta^2)(\alpha^2\sigma^2/2 - 6\alpha\beta\sigma + 2(1+5\beta^2))\sigma^{-3}.
\end{split}\end{align*}
The coefficient $c_1$ is the group velocity of the wave train, which corresponds to the speed at which the envelope of the wave train propagates through space. By choosing an appropriate co-moving frame, i.e. by taking $c = 2(\alpha - \beta)$ in~\eqref{newcoord}, we can factor out this transport. It depends on the sign of the diffusivity coefficient $c_2$ whether the curve $\lambda_+[\R]$ touches the origin as a left- or right-oriented parabola. Thus, the sign of $c_2$ determines the stability of the spectrum close to the origin. For all $(\alpha,\beta) \in \R^2$ satisfying $1+\alpha\beta < 0$ it holds $c_2 < 0$, implying the wave train is spectrally unstable. On the other hand, for $(\alpha,\beta) \in \R^2$ with $1+\alpha\beta > 0$ we have $c_2 > 0$ if and only if
\begin{align*}
\sigma > \sigma_{s} := \frac{2\left(1+\beta^2\right)}{1+\alpha\beta},
\end{align*}
or equivalently, using~\eqref{newcoord},
\begin{align*}
\zeta^{2} < \zeta_{s}^2, \qquad \zeta_s := \sqrt{\frac{1+\alpha \beta }{\left( 2\left( 1+\beta^{2}\right) +1+\alpha \beta \right) }}.
\end{align*}
Hence, as $|\zeta|$ is increased through $\zeta_{s}$ the wave-train solution $\Psi_\per$ undergoes a sideband instability. Although such a sideband instability occurs for all $(\alpha,\beta) \in \R^2$ with $1+\alpha\beta > 0$, the wave train is only {destabilized} through the sideband instability if and only if $(\alpha,\beta)$ is contained in some set $\A_s$, which is explicitly determined in~\cite{vH94}. There, one defines the function $r \colon (0,1) \to \R \cup \{\infty\}$ by $r(z) = z^{-1/2}$ for $z \in (0,1/3]$, by $r(z) = \infty$ for $z \in [3/4,1)$ and by the unique positive real root $r(z)$ to the algebraic equation
\begin{align*} r^4 z^2(4z-3) + r^2(5z^2 - 4z + 1) + 1 = 0,\end{align*}
for $z \in (1/3,3/4)$. Subsequently, one obtains
\begin{align*} \mathcal{A}_s & = \left\{(\alpha,\beta) \in \R^2 : \left(-1 < \alpha\beta < \beta^2\right) \vee (\beta = 0 \wedge \alpha \neq 0)\right.\\
&\qquad \qquad \qquad \qquad \qquad \qquad \qquad \qquad \left.\vee \left(0 < |\beta| \leq |\alpha| < r\left(\tfrac{\beta}{\alpha}\right)\right)\right\},\end{align*}
see also Figure~\ref{vanHarten}. Thus, for any $(\alpha,\beta) \in \mathcal{A}_s$ the wave train destabilizes through a sideband instability as $|\zeta|$ is increased through $\zeta_{s}$, which implies that for such $(\alpha,\beta)$ the Eckhaus boundary is given by $\zeta_{\bd} = \zeta_s$. At the Eckhaus boundary $|\zeta| = \zeta_{\bd}$, or equivalently at $\sigma = \sigma_s$, the expression for $c_3$ simplifies to
\begin{align}
c_{3,s} = 2\left(\alpha - \beta\right)\sigma_s^{-1}. \label{c3Eck}
\end{align}
So, if $(\alpha,\beta) \in \mathcal{A}_s$ with $\alpha \neq \beta$, we expect that the leading-order behavior of $\partial_t - \El$ is dispersive near the Eckhaus boundary $|\zeta| = \zeta_{\bd}$, which,  in the appropriate co-moving frame, will lead to a KdV equation as an amplitude equation for modulations of $\Psi_\per$ solving~\eqref{modeq2} -- see~\cite{vH94} and~\S\ref{sec3}. On the other hand, in the special case $\alpha = \beta$, $c_{3,s}$ vanishes and a Cahn-Hilliard equation occurs instead -- see again~\cite{vH94}.

\begin{figure}[h]
 \centering
  \begin{tikzpicture}[rotate=0,scale=1]
    \begin{axis}[
	xmin=-1.5, xmax=1.5,
	ymin=-0.6, ymax=0.5,
	axis lines=center,
	ticks=none,
    ]
      \addplot+[black,thick, domain=-2:2,samples=100,no marks,opacity=0.5] {-x^2*(x^2-1)};
      \node[below,black] at (axis cs:0.8,-0.02) {$ \mathcal{O}(\varepsilon) $};
           \node[below,black] at (axis cs:1.4,-0.02) {$ k $};
      \node[] at (axis cs:-0.35,0.4) {$\textrm{Re}(\lambda_+)$};
      \draw[] (axis cs:-0.05,0.25) -- (axis cs:0.05,0.25);
      \node[left=4pt] at (axis cs:0.05,0.26) {$ \mathcal{O}(\varepsilon^4) $};
    \end{axis}
  \end{tikzpicture}
 \caption{The spectral curve $ \lambda_+ $   in the parameter regime~\eqref{parreg}.}
 \label{sidebandfig}
\end{figure}
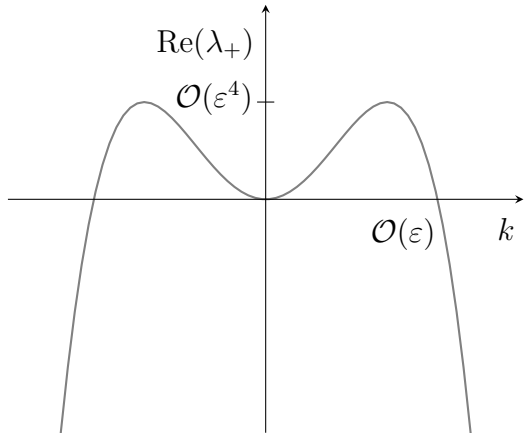

We are interested in proving the validity of the KdV equation as a modulation equation for marginally sideband unstable wave trains such that a small branch of unstable spectrum attached to the origin lies the right-half plane -- see Figure~\ref{sidebandfig} and Remark~\ref{remm}.

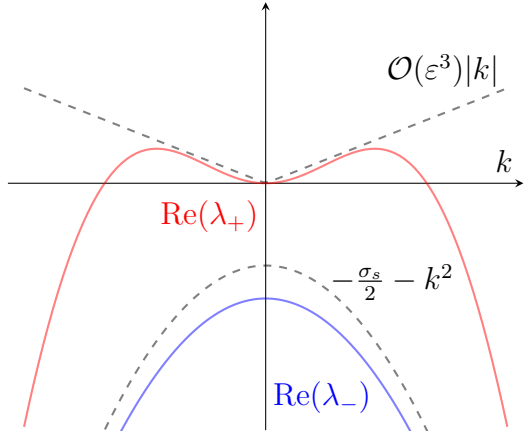
\begin{figure}[h]
 \centering
  \begin{tikzpicture}[rotate=0,scale=1]
    \begin{axis}[
	xmin=-1.6, xmax=1.6,
	ymin=-1.5, ymax=1.1,
	axis lines=center,
	ticks=none,
	xlabel={$k$},
    ]
      \addplot+[black,thick,dashed, domain=-1.5:1.5,samples=100,no marks,opacity=0.5] {0.385*abs(x)};
      \addplot+[red,thick, domain=-1.5:1.5,samples=100,no marks,opacity=0.5] {-x^2*(x^2-1)/(1+0.4*x^2)};
      \node[below,black] at (axis cs:1.1,0.85) {$\mathcal{O}(\varepsilon^3) |k|$};
      \addplot+[blue,thick, domain=-1.5:1.5,samples=100,no marks,opacity=0.5] {-x^2-0.7};
      \addplot+[black,thick,dashed, domain=-1.5:1.5,samples=100,no marks,opacity=0.5] {-x^2-0.5};
      \node[red] at (axis cs:-0.35,-0.2) {\red $\textrm{Re}(\lambda_+)$};
      \node[left=4pt] at (axis cs:1.3,-0.6) {$-\frac{\sigma_s}{2} -  k^2$};
      \node[blue] at (axis cs:0.35,-1.3) {$\textrm{Re}(\lambda_-)$};
    \end{axis}
  \end{tikzpicture}
 \caption{Spectral estimates established in Lemma~\ref{specboundslem}.}
 \label{fig:spectrum-linear-part-hypothesis}
\end{figure}

Therefore, we set
\begin{align}(\alpha,\beta) \in \A_s, \quad \text{ with } \alpha \neq \beta, \quad \sigma = \sigma_s - \varepsilon^2 = \frac{2\left(1+\beta^2\right)}{1+\alpha\beta} - \varepsilon^2, \label{parreg}\end{align}
where $0 < \varepsilon \ll 1$ is a small parameter. Recalling $\sigma = \zeta^{-2}-1$ by~\eqref{newcoord}, $|\zeta|$ lies just above the Eckhaus boundary $\zeta_{\bd}$ in the parameter regime~\eqref{parreg}, i.e., it holds $0 < |\zeta| - \zeta_{\bd} = \mathcal{O}(\varepsilon^2)$.
\begin{remark}\label{remm}{\rm
For notational simplicity we restrict ourselves to the marginally sideband-unstable case $0 < |\zeta| - \zeta_{\bd}  $.
One readily observes that  our analysis works in the marginally sideband-stable case $0 \geq |\zeta| - \zeta_{\bd} = \mathcal{O}(\varepsilon^2)$, too.
}\end{remark}
By~\eqref{c3Eck} and~\eqref{parreg} we have
\begin{align*} c_3 = c_{3,s} + \mathcal{O}(\varepsilon^2) \neq 0,\end{align*}
yielding the desired dispersive dynamics on the linear level. Finally, the parameter regime~\eqref{parreg} leads to the following spectral bounds -- see also Figure~\ref{fig:spectrum-linear-part-hypothesis}.

\begin{lemma} \label{specboundslem}
Assume~\eqref{parreg} is satisfied. There exists an $\epsilon$-independent constant $ C > 0$ such that the spectral curves $\lambda_\pm \colon \R \to \C$ given by~\eqref{speccurv} enjoy the following bounds
\begin{align*}
\Re(\lambda_-(k)) \leq -\frac{\sigma_s}{2} - k^2, \qquad \Re(\lambda_+(k)) \leq C \epsilon^3|k|, \qquad k \in \R,
\end{align*}
provided $\epsilon > 0$ is sufficiently small.
\end{lemma}
\noindent
{\bf Proof.}
Since $\upsilon(k)$ is a principal square root, it has positive real part. Hence, the bound on $\lambda_-$ follows immediately from~\eqref{speccurv} and~\eqref{parreg}. The function $f(k,\sigma) = \Re(\lambda_+(k;\sigma))$ depends smoothly on $k$ and $\sigma$ at $(k,\sigma) = (0,\sigma_s)$ and by~\eqref{critdisp} and~\eqref{c3Eck} it holds
\begin{align*} \partial_k^2 f(0,\sigma_s) &= 0, \qquad \partial_k^j f(0,\sigma) = 0, \qquad j = 0,1,3, \\
c_{4,s}  := -\partial_k^4 f(0,\sigma_s) & = \frac{(1+\alpha \beta)J(\alpha,\beta)}{4(1+\beta^2)^2}> 0, \\
J(\alpha,\beta) & = 1+5 \beta^2 + \alpha^2 (1-3 \beta^2) + 4 \alpha \beta (\beta^2 -1) ,
\end{align*}
for all $\sigma > 0$, where we use that $ 1+\alpha \beta > 0 $ and
$ J(\alpha,\beta)> 0 $
for $ (\alpha,\beta) \in \A_s $.
Thus, by Taylor's Theorem there exists a constant $C > 0$ such that
\begin{align*}
\left|\partial_k^2 f(0,\sigma)\right|, \left|\partial_k^4 f(0,\sigma) + c_{4,s}\right| \leq C\epsilon^2,\end{align*}
for $\epsilon > 0$ sufficiently small. On the other hand, we find -- again by Taylor's Theorem -- an $\epsilon$-independent neighborhood $U \subset \R$ of $0$ such that
\begin{align*}
\left|f(k,\sigma) - \frac12 \partial_k^2 f(0,\sigma)k^2 - \frac{1}{4!}\partial_k^4 f(0,\sigma) k^4\right| \leq \frac{c_{4,s} k^4}{2 \cdot 4!}, \qquad k \in U,
\end{align*}
as long as $\epsilon > 0$ remains bounded. Combining the latter two inequalities yields
\begin{align}
f(k,\sigma) \leq -\frac{c_{4,s} k^4}{4 \cdot 4!} + 2C\epsilon^2 k^2 \leq \frac{32C\sqrt{C}}{3\sqrt{c_{4,s}}} \epsilon^3|k|, \qquad k \in U, \label{specineq1}
\end{align}
for $\epsilon > 0$ sufficiently small, where we used that for any $a,b > 0$, the lines $k \mapsto \pm 2b\sqrt{b}k/(3\sqrt{3a})$ are tangent to the quartic $k \mapsto -ak^4 + bk^2$ at $k = \pm\sqrt{b/(3a)}$. On the other hand, since we have $(\alpha,\beta) \in \A_s$, the wave train undergoes a sideband instability at $\sigma = \sigma_s$, i.e., $f(k;\sigma_s)$ is strictly negative for all $k \in \R\setminus\{0\}$ and touches the origin at $k = 0$ in a quartic tangency. Thus, since $U$ is an $\epsilon$-independent neighborhood of the origin, we have $f(k,\sigma) < 0$ for $k \in \R \setminus U$, provided $\epsilon > 0$ is sufficiently small. Combining the latter with~\eqref{specineq1} yields the result.
\qed

\section{Derivation of the KdV equation}
\label{sec3}

In this section, we formally show that the KdV equation describes the leading-order behavior of modulations of the wave-train solution~\eqref{modsol} in the parameter regime~\eqref{parreg}. The KdV equation is a long-wave approximation, i.e., in Fourier space the solutions are localized about $ k = 0 $. Thus, we assume~\eqref{parreg} and make the ansatz that small solutions $(\psi_{*,\epsilon},s_{*,\epsilon})(x,t)$ to the  modulation equation~\eqref{modeq3}  have the form
\begin{align}
\begin{split}
\psi_{*,\epsilon}(x,t) &= \varepsilon^2 A(\varepsilon x , \varepsilon^3 t ) , \qquad s_{*,\epsilon}(x,t) = \varepsilon^2 B(\varepsilon x , \varepsilon^3 t ),
\end{split} \label{KdVAnsatz}\end{align}
with small parameter $0 < \varepsilon \ll 1$.
We expect that $ B $ is slaved by $ A $ and so we refine the ansatz to
\begin{align}
B(\xi,\tau) = \nu_0 A(\xi,\tau) + \varepsilon \nu_1 \partial_{\xi} A(\xi,\tau) + \varepsilon^2 \nu_2 \partial_{\xi}^2 A(\xi,\tau) + \varepsilon^2 \nu_3 (A(\xi,\tau))^2,
\label{KdVAnsatz2}
\end{align}
with  coefficients $\nu_{0,1,2,3} \in \R$, which will be determined later. Our goal is to formally derive that $A$ satisfies a KdV equation.

Inserting the ansatz~\eqref{KdVAnsatz} into~\eqref{modeq3} leads to the system
\begin{align}
 - \varepsilon^3 c \partial_{\xi} A + \varepsilon^5 \partial_{\tau} A &=
\varepsilon^4 \partial_{\xi}^2 A + \varepsilon^5 \alpha \partial_{\xi}^3 B
-\beta \sigma \partial_{\xi}(  (2 \varepsilon^3 B + 2\varepsilon^5 B^2 +  \mathcal{O}(\varepsilon^7)  ) \nonumber \\ & + \alpha   \varepsilon^7 \partial_{\xi}((\partial_{\xi} B)^2)  - 2 \varepsilon^3 \alpha    \partial_{\xi} A-   \varepsilon^5 \alpha  \partial_{\xi}(A^2) \nonumber
\\ & +2 \varepsilon^4     \partial_{\xi}^2 B  +2 \varepsilon^6  \partial_{\xi}(A(\partial_{\xi} B)) , \label{formeq}
\\
- \varepsilon^3 c \partial_{\xi} B + \varepsilon^5 \partial_{\tau} B  &=  \varepsilon^4 \partial_{\xi}^2 B - \varepsilon^3 \alpha \partial_{\xi} A - \sigma (2 \varepsilon^2 B + 2\varepsilon^4 B^2 + \mathcal{O}(\varepsilon^6) ) \nonumber \\ &
+  \varepsilon^6 (\partial_{\xi} B)^2- 2 \varepsilon^2  A-   \varepsilon^4 A^2
-2 \varepsilon^3 \alpha     \partial_{\xi} B  -2 \varepsilon^5 \alpha   A(\partial_{\xi} B) . \nonumber
\end{align}

Equating terms at order $\varepsilon^2$, $\varepsilon^3$ and $\varepsilon^4$ in the $B$-equation in~\eqref{formeq} to zero, while assuming for the moment that $\sigma$ is independent of $\varepsilon$, yields
\begin{align}\begin{split}
-  2 \sigma \nu_0
- 2    = 0,
\\
-  2 \sigma \nu_1
+ c  \nu_0
-  \alpha
- 2 \alpha \nu_0 = 0,
\\
-  2 \sigma \nu_2
+  \nu_0
+ c   \nu_1
-2  \alpha      \nu_1 = 0,
\\
-  2 \sigma \nu_3
- 2 \sigma \nu_0^2
- 1=0.
\end{split} \label{eqBeq}\end{align}
Similarly, equating terms at order $\varepsilon^3$, $\varepsilon^4$ and $\varepsilon^5$ in the $A$-equation in~\eqref{formeq} to zero leads to
\begin{align}
 c
- 2 \beta  \sigma   \nu_0
- 2  \alpha  =0 , \label{eqAeq1}
\\
1
+2         \nu_0
- 2 \beta  \sigma \nu_1=0,
 \label{eqAeq}\end{align}
and to the KdV equation
\begin{equation}\label{eq:kdv}
 \partial_{\tau} A  = \widetilde{\gamma}_{lin} \partial_{\xi}^3 A  +  \widetilde{\gamma}_{non} \partial_{\xi}(A^2),
\end{equation}
with coefficients
\begin{align*}\begin{split}
\widetilde{\gamma}_{lin} & :=
 \alpha \nu_0
 - 2 \beta   \sigma  \nu_2
 +2    \nu_1, \\ 
\widetilde{\gamma}_{non} & := - (2\beta   \sigma  \nu_0^2
  +  \alpha
 +2 \beta   \sigma   \nu_3 ).
\end{split} \end{align*}
Our aim is to express these coefficients in terms of $\alpha,\beta$ and $\sigma$ by solving~\eqref{eqBeq} with respect to $\nu_{0},\ldots,\nu_{3}$.
First, we choose the co-moving frame in~\eqref{newcoord} such that~\eqref{eqAeq1} is satisfied, i.e., we take the velocity
\begin{align}
c = 2 \beta  \sigma   \nu_0 + 2  \alpha. \label{veloc}
\end{align}
With this choice of $c$, system~\eqref{eqBeq} uniquely determines the coefficients $\nu_{0},\ldots ,\nu_{3}$. Indeed, we find
\begin{align}\begin{split}
\nu_0 &= - \sigma^{-1},\\
2 \sigma  \nu_1 &= 2 \beta  \sigma^{-1} -  \alpha,\\
2 \sigma  \nu_2 &=  -\sigma^{-1} - 2 \beta \nu_1,\\
2 \sigma  \nu_3 &= - 2\sigma^{-1} - 1.
\end{split} \label{coeffB}\end{align}
Substituting $\nu_0 = -\sigma^{-1}$ into~\eqref{veloc}, we recover the group velocity
\begin{align} c = 2(\alpha - \beta), \label{veloc2}\end{align}
of the critical wave number $k = 0$, meaning we have switched to a co-moving frame in which the envelope of the wave train does not propagate -- see~\S\ref{sec2}. Finally, by our choice $\sigma = \sigma_s - \varepsilon^2$ in~\eqref{parreg}, equation~\eqref{eqAeq} is satisfied to leading order. Indeed, we find
\begin{align*}\begin{split}
1+ 2 \nu_0 - 2 \beta  \sigma  \nu_1 = 1 - 2  \sigma_{s}^{-1} - 2 \beta^2  \sigma_{s}^{-1} +  \alpha \beta + \mathcal{O}(\varepsilon^2) = \mathcal{O}(\varepsilon^2).
\end{split}\end{align*}
Using~\eqref{parreg} and~\eqref{coeffB}, we approximate the coefficient $\widetilde{\gamma}_{lin}$ of the linear term in the KdV equation~\eqref{eq:kdv} by
 \begin{align*}\begin{split}
\widetilde{\gamma}_{lin}& =   \alpha \nu_0
 - 2 \beta  \sigma \nu_2
 +2     \nu_1\\
 &= ( \beta - \alpha )\sigma^{-1} + 2 \left(\beta^2 +1\right) \nu_1\\
 &= 2(\beta - \alpha) \sigma_s^{-1} + \mathcal{O}(\varepsilon^2)\\
 &= \gamma_{lin} + \mathcal{O}(\varepsilon^2) \neq 0,
 \end{split}\end{align*}
where $\gamma_{lin} := -c_{3,s} = (\beta - \alpha)(1+\alpha\beta)/(1+\beta^2)$ is the leading-order $c_3$-coefficient in the expansion~\eqref{critdisp} at the Eckhaus boundary -- see~\eqref{c3Eck}. Similarly, we approximate the coefficient $\widetilde{\gamma}_{non}$ in the KdV equation~\eqref{eq:kdv}
 \begin{align*}\begin{split}
\widetilde{\gamma}_{non}& =  - (2\beta  \sigma  \nu_0^2
  +  \alpha
 +2 \beta \sigma   \nu_3 )
=  \gamma_{non} + \mathcal{O}(\varepsilon^2) \neq 0,
\end{split}\end{align*}
where $\gamma_{non} := \beta - \alpha$. We conclude that with the choice of coefficients~\eqref{coeffB} and velocity~\eqref{veloc2} the equations~\eqref{eqBeq} are satisfied and the equations~\eqref{eqAeq1} and~\eqref{eqAeq} are satisfied to leading order with $\mathcal{O}(\varepsilon^2)$ residual. Thus, taking $A$ as a solution to the KdV equation~\eqref{KdV}, we find that the ansatz in~\eqref{KdVAnsatz}-\eqref{KdVAnsatz2} formally solves the modulation  equation~\eqref{modeq3}  at least up to order $\mathcal{O}(\varepsilon^4)$.
In order to prove that the KdV equation~\eqref{KdV} makes correct predictions
about the dynamics of~\eqref{eq1}
this has to be improved
subsequently in~\S\ref{sec:imprAnsatz} by adding higher order terms to the approximation.
However, the construction of the improved approximation will be made in a more adequate
coordinate system.


\section{The functional analytic set-up}
\label{setup}

In order to  state our main result we  introduce a number of function spaces and notations.
By $ \langle \cdot, \cdot \rangle $ we denote the Euclidean inner product and by
$|\cdot|$  the associated Euclidean norm in $ \R^d $.
The Fourier transform is denoted by
\begin{align*}
\F(u)(k) = \hat{u}(k) = \frac{1}{\sqrt{2\pi}} \int_\R e^{-ikx} u(x)dx.
\end{align*}
For $ m \geq 0 $ we   define the Sobolev spaces
$$
H^m = \{u \in L^2(\R) : (1+|\cdot|^{2})^{\frac{m}{2}}\hat{u} \in L^2(\R)\},
$$
endowed with the inner product
\begin{align*} \langle u, v\rangle_{H^m} =
 \langle \widehat{u}, \widehat{v}\rangle_{L^2_m} =
\int_{\R} \left(1+|k|^{2}\right)^m \langle \hat{u}(k),\hat{v}(k)\rangle \; dk.\end{align*}
For any $m \in \N$, the induced norm
is equivalent to the usual $ H^m $-norm.
Finally, for $ m \geq 0 $ we introduce
$$
W_{m} :=  \left\{u : u = \F^{-1}(\widehat{u}), \widehat{u} \in L^1(\R), \|u\|_{W_m} = \int_\R (1+|k|^{m})|\widehat{u}(k)|\; dk < \infty\right\}.
$$
By Sobolev's embedding theorem the space
$H^{m+\delta}(\R)$ is continuously embedded into $W_{m}$ for each $\delta > 1/2$.
Moreover, every $ u \in W_{m} $ is $ \lfloor m \rfloor $-times continuously differentiable with finite $  C^{\lfloor m \rfloor}_b(\R) $-norm.

In the parameter regime~\eqref{parreg}, the wave train is marginally sideband-unstable, see Figure~\ref{sidebandfig}, leading to positive growth rates of the semigroup associated to the linearization. To account for these growth rates, we work in the space
$$
H^\infty_{\mu,s} = \{u \in L^2(\R) : e^{\mu|\cdot|}(1+|\cdot|^2)^{\frac{s}{2}}\hat{u} \in L^2(\R)\},
$$
endowed with the norm
\begin{align*}
\|u\|_{H^\infty_{\mu,s}} = \left(\int_{\R} |\hat{u}(k)|_2^2 e^{2\mu|k|} (1+|k|^2)^s dk\right)^{\frac{1}{2}},
\end{align*}
where $\mu \geq  0$ and  $s  \geq 0$. Functions $ u \in H^\infty_{\mu,0}$  can be extended to functions that are analytic on the strip $\{z \in \C \colon |\Im(z)| < \mu\}$.
In the following we use the abbreviation $ H^\infty_{\mu} = H^\infty_{\mu,0} $.
It is readily seen that for any $\mu_1  > \mu_2 \geq 0 $ and any $m \geq 0$ we have the
 continuous embedding
$H^\infty_{\mu_1,0} \subset H^\infty_{\mu_2,m} $.

Similarly, we define the spaces $ W_{\mu,m} $.

In our notations of the spaces and norms we do not distinguish between
scalar and vector-valued functions.

\section{Main results} \label{mainresults}

In~\S\ref{sec3} we formally derived a  KdV approximation for small long-wave modulations of wave-train solutions to the GL equation~\eqref{eq1} in the parameter regime~\eqref{parreg}. In this section, we state rigorous approximation results showing that the KdV equation makes indeed correct predictions about the dynamics of the modulated wave train on a non-trivial time scale.

In order to state our main approximation results, we assume~\eqref{parreg} and switch to a comoving frame~\eqref{newcoord}, where the velocity $c$ is given by~\eqref{veloc2}. We consider an $L^2$-solution $A$ to the KdV equation~\eqref{KdV}, which is analytic on a strip in the complex plane. We emphasize that such solutions exist locally in time, see Theorem~\ref{linkin1}. The associated long-wave solution~\eqref{KdVAnsatz}-\eqref{KdVAnsatz2} with coefficients~\eqref{coeffB} provides a non-trivial approximation of an $\mathcal{O}(\varepsilon^2)$-solution to the modulation equation~\eqref{modeq3}  on the long $\mathcal{O}(1/\epsilon^{3})$-time-scale.

\begin{theorem} \label{mainresult1}
Let $m, \mu_{\A}, \tau_0 > 0$ and $A \in C([0,\tau_0],H^\infty_{\mu_{\A}})$ be a solution to the KdV equation~\eqref{KdV}. Then there exists $C,\tau_1, \varepsilon_0 > 0$ such that for all $\epsilon \in (0,\epsilon_0)$ a solution $V(x,t) = (\psi,s)(x,t)$ to the modulation equation~\eqref{modeq3} exists with
\begin{align*} \sup_{0 \leq \epsilon^3 t \leq \tau_1} \left\|V(\cdot,t) - V_{app}^{*,\epsilon}(\cdot,t)\right\|_{H^m} \leq C\epsilon^{5/2}, \end{align*}
where $V_{app}^{*,\epsilon}(x,t) = (\psi_{app}^{*,\epsilon},s_{app}^{*,\epsilon})(x,t)$ is defined by~\eqref{KdVAnsatz}-\eqref{KdVAnsatz2} and~\eqref{coeffB}. In particular, it holds
\begin{align*} 
  \sup_{0 \leq \epsilon^3 t \leq \tau_1} \ \sup_{x \in \R} |V(x,t) - V_{app}^{*,\epsilon}(x,t)| \leq C\epsilon^{3}. 
\end{align*}
\end{theorem}

The error bound in Theorem~\ref{mainresult1} can be improved by adding higher-order terms
to the approximation ansatz~\eqref{KdVAnsatz} -- see~\S\ref{sec:imprAnsatz}. This leads to the following statement.

\begin{theorem} \label{mainresult2}
Let $m,\mu_{\A}, \tau_0 > 0$, $\kappa > 3 $, and $A \in C([0,\tau_0],H^\infty_{\mu_{\A}})$ be a solution to the KdV equation~\eqref{KdV}. Then there exists $C,\tau_1, \varepsilon_0 > 0$ such that for all $\epsilon \in (0,\epsilon_0)$ there exists an approximation $V_{app}^{\kappa,\epsilon} \colon \R \times [0, \tau_1 /\epsilon^{3}] \to \R^2$ satisfying
\begin{align} \sup_{0 \leq \epsilon^3 t \leq \tau_1} \|V_{app}^{\kappa,\epsilon}(\cdot,t) - V_{app}^{*,\epsilon}(\cdot,t)\|_{H^m} &\leq C\epsilon^{5/2}, \label{impresolv}\end{align}
where $V_{app}^{*,\epsilon}(x,t) = (\psi_{app}^{*,\epsilon},s_{app}^{*,\epsilon})(x,t)$ is defined by~\eqref{KdVAnsatz}-\eqref{KdVAnsatz2} and~\eqref{coeffB} and there is a solution $V(x,t) = (\psi,s)(x,t)$ to the modulation  equation~\eqref{modeq3}
satisfying
\begin{align} \sup_{0 \leq \epsilon^3 t \leq \tau_1} \left\|V(\cdot,t) - V_{app}^{\kappa,\epsilon}(\cdot,t)\right\|_{H^m} &\leq C\epsilon^{\kappa}. \label{Hmbound}\end{align}
\end{theorem}
\begin{corollary} \label{mainresult2coro}
Under the assumptions of Theorem~\ref{mainresult2}
 it holds
\begin{align}
\sup_{0 \leq \epsilon^3 t \leq \tau_1} \ \sup_{x \in \R} |V(x,t) - V_{app}^{\kappa,\epsilon}(x,t)| & \leq C\epsilon^{\kappa}, \label{pointwise21}\\
\sup_{0 \leq \epsilon^3 t \leq \tau_1} \ \sup_{x \in \R} |V_{app}^{\kappa,\epsilon}(x,t) - V_{app}^{*,\epsilon}(x,t)| & \leq C\epsilon^3. \label{pointwise22}
\end{align}
\end{corollary}
\noindent {\bf Proof.}
The $L^\infty$-bound~\eqref{pointwise21} follows from~\eqref{Hmbound}, since $H^m$ is continuously embedded into $L^\infty(\R)$ for $ m > 1/2 $. The bound~\eqref{pointwise22} follows by construction of the improved approximation ansatz in~\S\ref{sec:imprAnsatz}-\S\ref{sec:imprAnsatz2} and
the calculations in Remark~\ref{compareAnsatz}, noting that for $0 \leq \tau \leq \tau_1$ the difference $V_{app}^{\kappa,\epsilon}(\cdot/\epsilon,\tau/\epsilon^3) - V_{app}^{*,\epsilon}(\cdot/\epsilon,\tau/\epsilon^3)$ is of order $\ord(\epsilon^3)$ in $L^\infty(\R)$ for some $\mu_\A > 0$ -- see Theorem~\ref{th66} -- and the $L^\infty$-norm is invariant under rescaling.
\qed

Since Theorem~\ref{mainresult1} is a direct consequence of Theorem~\ref{mainresult2}
and Corollary~\ref{mainresult2coro}, it remains to prove Theorem~\ref{mainresult2} -- see~\S\ref{secproof}.

\begin{remark}{\rm
Our approximation result is not optimal in the sense that
in general $ \tau_1  < \tau_0$. However,  we still obtain the natural and non-trivial $ \mathcal{O}(1/\varepsilon^3) $-time scale for the approximation time.
 }
\end{remark}

\begin{remark}{\rm
Above approximation results should not be taken for granted. There are counterexamples
that  formally derived amplitude  equations make wrong predictions, cf.~\cite{Schn96MN} and~\S\ref{seccounterexample}.
}\end{remark}

\begin{remark}\label{rem65}{\rm
To construct the improved KdV approximation $V_{app}^{\kappa,\epsilon} $ from the original KdV approximation $V_{app}^{*,\epsilon} $, we have to solve, additional to the KdV equation, a number of inhomogeneous linearized
KdV equations, cf.~\S\ref{sec:imprAnsatz}.
Besides in the proof of Theorem~\ref{mainresult1} and Theorem~\ref{mainresult2}
the improved approximation $V_{app}^{\kappa,\epsilon} $ is utilized
to transfer the approximation of solutions to the modulation equation~\eqref{modeq3} in Theorem~\ref{mainresult2}  to the original modulated solution~\eqref{modsol} to the complex Ginzburg-Landau equation. This is discussed in detail in~\S\ref{secbackU}.
}
\end{remark}

%

\section{Proof of  Theorem~\ref{mainresult2}} \label{secproof}

Without loss of generality we assume $m \geq 2$. The result for smaller $m$-values is an immediate consequence.  The bound~\eqref{impresolv} follows  by construction of the improved approximation ansatz in~\S\ref{sec:imprAnsatz}-\S\ref{sec:imprAnsatz2}.
The error bound~\eqref{Hmbound} will be proven in~\S\ref{sec:resest}-\S\ref{sec:error}.
Before we do so, we outline below some more details, especially we motivate the change of variables made in~\S\ref{sec:diagonal}. A reformulation of Theorem~\ref{mainresult2}
in the new variables
can be found in~\S\ref{secmainresult7} and the transfer of the result back to the original variables in~\S\ref{sec78}.

\subsection{Outline} \label{sec:outline}

As already said, it is a non-trivial task to bound solutions $V = (\psi,s)$ to the modulation  equation~\eqref{modeq3}
of order
$ \mathcal{O}(\varepsilon^2) $ on the long $ \mathcal{O}(1/\varepsilon^{3}) $-time interval or, equivalently, to estimate the error $ R_V $ defined  by
\begin{align*}
V = V_{app}^{\kappa,\epsilon} + \epsilon^{\kappa} R_V
\end{align*}
in $H^m$ by an $\epsilon$-independent bound on the long $\ord(1/\epsilon^{3})$-time scale.

Our approach to tackle this problem is to establish an $\epsilon$-independent constant $C>0$ and a differential inequality of the form
\begin{align} \E'(t) \leq C \epsilon^3 (\E(t) + 1), \label{diffineq}\end{align}
for a suitably chosen energy $\E(t)$, which allows to estimate  $\|R_V(t)\|_{H^m}$.  Applying Gr\"onwall's inequality to~\eqref{diffineq} leads then to the desired $H^m$-bound on $R_V(t)$ on an $\ord(1/\epsilon^{3})$-time scale. We outline below that one has to overcome a number of problems.

Since $V$ solves~\eqref{modeq3}, the error $R_V$ satisfies
\begin{align} 
  \partial_t R_V= L R_V + G (V_{app}^{\kappa,\epsilon},R_V) - \epsilon^{-\kappa} \Res_V(V_{app}^{\kappa,\epsilon}), \label{erroreq}
\end{align}
where  $G$ contains linear and nonlinear terms with respect to $ R_V$ and is given by
\begin{align*} 
  G (V_{app}^{\kappa,\epsilon},R_V) = \epsilon^{-\kappa} \left(N(\epsilon^\kappa R_V + V_{app}^{\kappa,\epsilon}) - N(V_{app}^{\kappa,\epsilon})\right),
\end{align*}
and where the {residual} $\Res_V(V)$ is defined by
\begin{align}
\Res_V(V) = \partial_t V - LV - N(V).
\label{resform}
\end{align}

Thus, in order to obtain a differential inequality of the form~\eqref{diffineq}, we require an $\ord(\epsilon^{3+\kappa})$-bound on the residual. The nonlinear terms with respect to $R_V$ in $G$ are of order $\ord(\epsilon^\kappa)$, so that we require $\kappa > 3$. The major difficulty comes from the linear terms with respect to $R_V$ in $G$. Because the KdV approximation $V_{app}^{\kappa,\epsilon}$ is of order $\ord(\epsilon^2)$, these linear terms are proportional to $\ord(\epsilon^{2})$. In order to extract the properties of our system, which nevertheless allows to establish an inequality of the form~\eqref{diffineq} for a suitable chosen energy, we have to perform additional changes of variables.

To gain the same regularity in both components of $R_V$, we first modify the local wave number $\psi = \partial_x \phi$, which was introduced in~\S\ref{sec1b}. We replace $\psi = \partial_x \phi$ in~\eqref{modeq3} by $\chi = \theta(\phi)$, which is defined through its Fourier symbol $\hat{\theta} \colon \R \to \R$ given by $\hat{\theta}(k) = ik \min\{1,|k|^{-1}\}$ and is called local pseudo-derivative in the following -- see Figure~\ref{fig:pseudo-derivative}. After applying this transform, the linearity $L$ in the error equation~\eqref{erroreq} is replaced by a linearity that has the same spectrum as the linear operator $\El$ in system~\eqref{modeq2}. Thus, by the calculations in~\S\ref{sec2}, it has unstable spectrum close to the origin in the parameter regime~\eqref{parreg} -- see also Figure~\ref{sidebandfig} -- but also an exponentially damped part. This exponentially damped part allows us to control $\ord(\epsilon^2)$-terms in $G$ on the long $\ord(1/\epsilon^{3})$-time scale. The weakly unstable part is of no help in this respect: one has to use that the corresponding part in $G$ has a local pseudo-derivative in front -- see also~\eqref{nonl1}. In order to exploit this fact we introduce a time-dependent exponential weight, which damps at all wave numbers of the weakly unstable part except for the wave number $k = 0$. This is sufficient to control the corresponding $\ord(\epsilon^2)$-terms in $G$ that have a pseudo-derivative in front, since they vanish at $k = 0$. Such time-dependent exponential weights are also employed in the functional analytic versions of the Cauchy-Kowalevski theorem, cf.~\cite{Ovs76}.

\subsection{The change of variables}
\label{sec:diagonal}

As outlined in~\S\ref{sec:outline}, we perform some changes of variables to the modulation equation~\eqref{modeq3}.

\begin{figure}[h]
 \centering
  \begin{tikzpicture}[rotate=0,scale=1]
    \begin{axis}[
	xmin=-2.3, xmax=2.3,
	ymin=-1.5, ymax=1.5,
	axis lines=center,
	ylabel={$-i \widehat \theta( k)$},
	xlabel={$k$},
    ]
      \addplot+[black,thick, domain=-1:1,samples=10,no marks,opacity=0.5] {x};
      \addplot+[black,thick, domain=-2:-1,samples=10,no marks,opacity=0.5] {-1};
      \addplot+[black,thick, domain=1:2,samples=10,no marks,opacity=0.5] {1};
    \end{axis}
  \end{tikzpicture}
 \caption{The symbol of the local pseudo-derivative.}
 \label{fig:pseudo-derivative}
\end{figure}
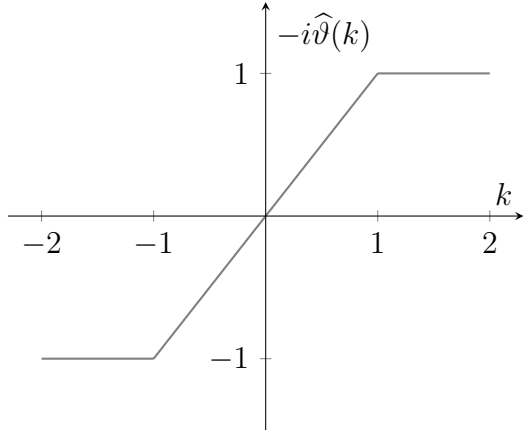

First, we replace $\psi = \partial_x \phi$ in~\eqref{modeq3} by the pseudo-derivative $\chi = \theta(\phi)$, which we define through its Fourier symbol $\hat{\theta}(k) = ik \min\{1,|k|^{-1}\}$, see Figure~\ref{fig:pseudo-derivative}.

This amounts to a transform $S_{\vartheta}$, which is defined through its Fourier symbol $
\widehat{S_{\vartheta}} $ given by
\begin{align*}
  \widehat{S_{\vartheta}}(k) = \begin{pmatrix} \min\{1,|k|^{-1}\} & 0 \\ 0 & 1\end{pmatrix}.
\end{align*}
We apply $S_{\vartheta}$ to~\eqref{modeq3} and find that $Y := S_{\vartheta}(V) =: (\chi,s)$ satisfies the equation
\begin{align} \partial_t Y = L_Y Y + N_Y(Y), \label{modeq4a}\end{align}
with linearity
\begin{align*}
L_Y = \begin{pmatrix}
  \partial_x^2 + (c-2\alpha) \partial_x & \theta \left(- 2\beta \sigma + 2\partial_x + \alpha\partial_x^2 \right)\\
  \theta^{-1}\left(-2 \partial_x -\alpha \partial_x^2\right) & - 2 \sigma + \partial_x^2 + (c-2\alpha) \partial_x
  \end{pmatrix},
\end{align*}
and nonlinearity
\begin{align}
\begin{split}
N_Y(Y) &= \begin{pmatrix}
	  \theta \left(- \beta \sigma h(s) + 2(\partial_x s)(\partial_x \theta^{-1} \chi) + \alpha (\partial_x s)^2 - \alpha (\partial_x \theta^{-1} \chi)^2 \right)\\
	  - \sigma h(s) + (\partial_x s)^2 - (\partial_x \theta^{-1} \chi)^2 - 2\alpha \partial_x s(\partial_x \theta^{-1} \chi)	
  \end{pmatrix}\\
 &=: \begin{pmatrix} \theta N_{Y1}(Y)\\ N_{Y2}(Y)\end{pmatrix},
 \end{split} \label{nonlN1}
\end{align}
where we recall $h(s) = e^{2s}-1-2s$. We note that $L_Y$ is a differential operator on $L^2(\R)$, whose highest proper derivatives are of second order. In addition, the nonlinearity $N_Y$ contains besides local pseudo-derivatives only first order proper derivatives. By construction, the first component of $ N_Y $  contains a local pseudo-derivative in front.

Since for fixed $k \in \R$ the Fourier symbols $\hat{L_Y}(k)$ and $\hat{\El}(k)$ of the operators $L_Y$ and $\El$ (note that $\El$ was defined as the linear part of system~\eqref{modeq2}) are similar matrices, the spectra of $L_Y$ and $\El$ coincide. Thus, using the spectral calculations in~\S\ref{sec2}, one observes that $L_Y$ has unstable spectrum close to the origin in the parameter regime~\eqref{parreg} but also an exponentially damped part -- see Figure~\ref{sidebandfig}.

The nonlinear
terms corresponding to
the exponentially damped part can be controlled on
the long $\ord(1/\epsilon^{3})$-time scale.
For the nonlinear
terms corresponding to
the weakly unstable part we exploit that they
have a local pseudo-derivative in front.

\begin{figure}[h]
 \centering
  \begin{tikzpicture}[rotate=0,scale=1]
    \begin{axis}[
	xmin=-1.5, xmax=1.5,
	ymin=-1.5, ymax=0.5,
	axis lines=center,
	ticks=none,
	xlabel={$ k $},
    ]
      \addplot+[red,thick, domain=-2:2,samples=250,no marks] {-x^2*(x^2-1)-0.7*abs(x)};
      \addplot+[blue,thick, domain=-2:2,samples=250,no marks] {-x^2-0.7-0.7*abs(x)};
      \node[red] at (axis cs:-0.35,-0.35) {$\textrm{Re}(\lambda_{2,+})$};
      \node[blue,anchor=west] at (axis cs:0.15,-0.8) {$\textrm{Re}(\lambda_{2,-})$};
    \end{axis}
  \end{tikzpicture}
 \caption{The transformed spectral curves  using the estimates from  Lemma~\ref{specboundslem}.}
 \label{fig:dampedspectrum}
\end{figure}
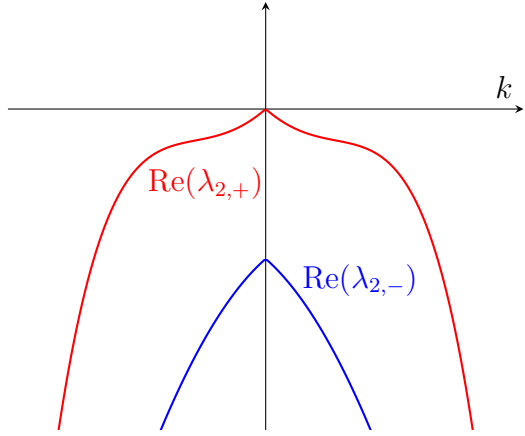

We introduce a time-dependent  exponential weight, which damps at all wave numbers of the weakly unstable part of $L_Y$ except for the wave number $k = 0$. We choose a $ \mu_* \in (0,\mu_{\A}) $, with
$\mu_{\A} > 0$ as in Theorem~\ref{mainresult2}.
Then we define the transform $S_{\omega}(t) \colon H^\infty_{\mu(t),m} \to H^m$ for $t\in[0,\mu_*/( \eta \epsilon^{3})]$,
with $ \mu(t) := (\mu_* - \eta \epsilon^3 t)/\epsilon $,
 through its Fourier symbol $\hat{S_{\omega}}(t,k) = \mathrm{diag}(\hat{\omega}(t,k),\hat{\omega}(t,k))
$, where  $\hat{\omega}(t,k) = e^{\mu(t)|k|}$ and   $\eta > 0$ is a $k$-, $\epsilon$- and $t$-independent constant yet to be defined. Applying $S_{\omega}(t)$ to~\eqref{modeq4a} we find that $\Y(\cdot,t) := S_{\omega}(t)Y(\cdot,t)$ satisfies
\begin{align} \partial_t \Y = L_\Y \Y + N_\Y(t,\Y), \label{modeq4b}\end{align}
with linearity
\begin{align*}
L_\Y = L_Y - \epsilon^2 \eta |k|_{op},
\end{align*}
where the pseudo-differential $|k|_{op} $ acts in Fourier space through multiplication by $|k|$. The nonlinearity in~\eqref{modeq4b} is given by
\begin{align*}
N_\Y(t,\Y) = S_{\omega}(t) N_Y(S_{\omega}(t)^{-1}(\Y)).
\end{align*}
The spectrum of $L_\Y$ is the union $\lambda_{2,+}[\R] \cup \lambda_{2,-}[\R]$, where $\lambda_{2,\pm}(k) := \lambda_\pm(k) - \epsilon^2 \eta |k|$ and $\lambda_\pm(k)$ is as in~\S\ref{sec2}. Lemma~\ref{specboundslem} provides the bounds
\begin{align}
\Re(\lambda_{2,-}(k)) \leq -\frac{\sigma_s}{2} - k^2, \qquad \Re(\lambda_{2,+}(k)) \leq -\epsilon^2\frac{\eta}{2} |k|, \qquad k \in \R, \label{eigvbounds}
\end{align}
provided $\epsilon > 0$ is sufficiently small. Thus, $L_\Y$ is damped for all wave numbers except for the wave number $k = 0$ -- see also Figure~\ref{fig:dampedspectrum}.
\begin{remark}{\rm
Note that the decay with $- \epsilon^2 \eta |k|$ is exploited to control the
$ \mathcal{O}(\varepsilon^2) $-terms of $ G $ in the error equation~\eqref{erroreq}. For making $\Re(\lambda_{2,+}(k)) $
negative a decay proportional to $ \epsilon^3  |k| $ would suffice.
}
\end{remark}

That the  nonlinear terms corresponding to the neutral mode really have a local pseudo-derivative in front, can be seen after diagonalizing the operator $L_\Y$ in Fourier space. We have the spectral decomposition
\begin{align*} \hat{L_\Y}(k) = \hat{S_{diag}}(k)^{-1}\hat{L_\cZ}(k)\hat{S_{diag}}(k),\end{align*}
for $k \in \R$ with
\begin{align}
\hat{L_\cZ}(k) = \mathrm{diag}\left(\lambda_{2,+}(k),\lambda_{2,-}(k)\right), \qquad k \in \R.\label{Fsymb}
\end{align}
and explicit representations for $ \hat{S_{diag}}(k) $ and
$ \hat{S_{diag}}(k)^{-1} $ in Appendix~\ref{secsdiag}.
In the following we use that  $ \hat{S_{diag}}(k) $ is of the form
$$
 \hat{S_{diag}}(k) =  \begin{pmatrix}
  s_{11}(k) &  \hat{\theta}(k) s_{12}(k)\\
  s_{21}(k) & s_{22}(k)
 \end{pmatrix},
$$
with smooth and bounded coefficients $ s_{ij}(k) $ and that
$ S_{diag} \colon H^m \to H^m$ given by $ S_{diag}= \F^{-1}[\hat{S_{diag}} \F(\cdot)]$ is an  isomorphism, cf. Lemma~\ref{applem-1}. Now, $\cZ := S_{diag}(\Y)$ satisfies a system of the form
\begin{align}
\partial_t \cZ = L_\cZ \cZ + N_\cZ(t,\cZ), \label{modeq4c}
\end{align}
with linearity $L_\cZ$ having the Fourier symbol $ \hat{L_\cZ} $.

So, $L_\cZ$ has, besides the wave number $k=0$, no unstable spectrum -- see~\eqref{eigvbounds}. The nonlinearity in~\eqref{modeq4c} is given by
\begin{align*}
\begin{split}
N_\cZ(t,\cZ) &= {S_{diag}}{N_\Y}(t,S_{diag}^{-1} \cZ)
\\ & =
S_{\omega}(t)  {S_{diag}}{N_Y}(S_{\omega}^{-1}(t) S_{diag}^{-1} \cZ)
\end{split}
\end{align*}
and has a local pseudo-derivative in front of the first component due to
\begin{equation} \label{psde}
\widehat{S_{diag}N_Y}  =  \begin{pmatrix}
  s_{11} &  \hat{\theta} s_{12}\\
  s_{21} & s_{22}
 \end{pmatrix}  \begin{pmatrix} \hat \theta \hat N_{Y1}\\ \hat N_{Y2}\end{pmatrix} 
 =  \begin{pmatrix}  \hat \theta \big(s_{11} \hat N_{11} +  s_{12} \hat N_{12} )\\
 s_{21} \hat \theta \hat N_{11} +  s_{22} \hat N_{12}
 \end{pmatrix},
\end{equation}
where we used~\eqref{nonlN1} and suppressed the $k$-dependency. Writing $\cZ = (\cZ_1,\cZ_2)$, we observe that the $\cZ_2$-equation is linearly exponentially damped, whereas the $\cZ_1$-equation is exponentially damped everywhere except for $k = 0$. Therefore, as outlined in~\S\ref{sec:outline}, it is crucial to use the conservation law structure of the critical modes, i.e.,  that the first term in the nonlinearity $N_{\cZ}$ has the local pseudo-derivative $\theta$ in front to compensate for the fact that $\lambda_{2,+}(0) = 0$ at wave number $k = 0$ of $L_{\cZ}$.
By Lemma~\ref{applem3} the nonlinear mapping $N_{\cZ}$ is smooth
from $ H^{s+1} $ to $ H^{s} $ for $ s \geq 2 $.

\begin{remark} \label{remresz} {\rm The diagonalization operator $ S_{diag} $ and the smoothing operator
$ S_{\omega}^{-1}(t) $ commute. For deriving higher order approximations in~\S\ref{sec:imprAnsatz}, see also Remark~\ref{rem65}, it turns out to be
advantageous to introduce
$Z := S_{diag}(Y)$, which  satisfies an equation of the form
\begin{align}
\partial_t Z = L_Z Z + N_Z(Z),
\label{modeq4aneu}
\end{align}
where $L_Z$  is a diagonal pseudo-differential operator  with Fourier symbol $\mathrm{diag}(\lambda_{+}(k),\lambda_{-}(k))$, where $ \lambda_{\pm} $ are as
in~\S\ref{sec2}. In addition,  the first component of
$ N_Z $  contains a local pseudo-derivative in front,
cf.~\eqref{psde}. We refer to Figure~\ref{figvariables} for an overview of all transformations.}
\end{remark}

\begin{figure}[htbp] 
   \centering
  \setlength{\unitlength}{0.6cm}
\begin{picture}(10, 6)(-1, -1)
\put(0,0){$ V $}
\put(5,0){$ Y $}\put(10,0){$ Z $} \put(5,4){$ \Y $}\put(10,4){$ \cZ $}
\put(1,0.3){\vector(3,0){3.3}}
\put(6,0.3){\vector(3,0){3.3}}
\put(6,4.3){\vector(3,0){3.3}}
\put(5.3,1){\vector(0,1){2.3}}
\put(10.3,1){\vector(0,1){2.3}}
\put(2,1){$ S_{\theta} $}
\put(7,1){$ S_{diag} $}
\put(7,5){$ S_{diag} $}
\put(3.4,2.2){$ S_{\omega}(t) $}
\put(10.8,2.2){$ S_{\omega}(t) $}
\end{picture}
   \caption{Summary of all transformations. In the equations for the variables in roman
 no smoothing occurs. In the equations for
   the calligraphic  variables a smoothing occurs. }
   \label{figvariables}
\end{figure}

\begin{remark}{\rm
The introduction of the time-dependent exponential weight  leads to  decay rates of the associated semigroups, which we do not have for  $  e^{t L_Z}: H^m \to H^m $.  We obtain now an estimate
\begin{align*}
 \| e^{t L_\cZ} \theta \|_{H^m \to H^m} \leq C(\epsilon,\eta) t^{-1}, \qquad t > 0,
\end{align*}
Thus, applying $ e^{t L_\cZ} $ to the first component of $ N_\cZ $ gives a polynomial decay rate. The latter could possibly be used to obtain the required estimates on the long $\ord(1/\epsilon^{3})$-time scale by using the variation of constant formula and optimal regularity theory. Instead of doing so we will work with energy estimates.}
\end{remark}

\begin{remark}{\rm
 There is a price we pay for the improved linear damping properties in the $\cZ$-equation~\eqref{modeq4c}. Starting from the residual $ \Res_V(t) $, defined in~\eqref{resform}, we introduce the residuals $\Res_{\cZ}(t) = S_{diag} S_\omega(t) S_\theta[\Res_V(t)]$ and $\Res_{Z} = S_{diag} S_\theta[\Res_V]$. In order to have the desired $\ord(\epsilon^{\kappa + 3})$-residual $\Res_{\cZ}$ in $H^m$ -- see~\S\ref{sec:outline} -- we require the residuals $ \Res_V $  and $ \Res_Z $ to be of order $\ord(\epsilon^{\kappa+3})$ in $H^\infty_{\mu_*/\epsilon,m}$ on  an $\ord(1/\epsilon^{3})$-time scale.}
\end{remark}

\subsection{The approximation result for the $ \cZ $-system}

\label{secmainresult7}

In order to prove Theorem~\ref{mainresult2} we first formulate the associated
approximation result for the $ \cZ $-system.
\begin{theorem} \label{mainresult7}
Let $m,\mu_{\A}, \tau_0 > 0$, $\kappa > 3 $, and $A \in C([0,\tau_0],H^\infty_{\mu_{\A}})$ be a solution to the KdV equation~\eqref{KdV}. Then there exists $C,\tau_1, \varepsilon_0 > 0$ such that for all $\epsilon \in (0,\epsilon_0)$ there exists an approximation $\cZ_{app}^{\kappa,\epsilon} \colon \R \times [0, \tau_1 /\epsilon^{3}] \to \R^2$ satisfying
\begin{align} \sup_{0 \leq \epsilon^3 t \leq \tau_1} \|\cZ_{app}^{\kappa,\epsilon}(\cdot,t) - S_{\omega(t)} S_{diag} S_{\theta} V_{app}^{*,\epsilon}(\cdot,t)\|_{H^m} &\leq C\epsilon^{5/2}, \label{impresolvz}\end{align}
where $V_{app}^{*,\epsilon}(x,t) = (\psi_{app}^{*,\epsilon},s_{app}^{*,\epsilon})(x,t)$ is defined by~\eqref{KdVAnsatz}-\eqref{KdVAnsatz2} and~\eqref{coeffB}, such that there is a  solution $\cZ $ to the modulation  equation~\eqref{modeq4c}
satisfying
\begin{align} \sup_{0 \leq \epsilon^3 t \leq \tau_1} \left\|\cZ(\cdot,t) - \cZ_{app}^{\kappa,\epsilon}(\cdot,t)\right\|_{H^m} &\leq C\epsilon^{\kappa}. \label{Hmboundz}\end{align}
\end{theorem}

\subsection{The improved approximation ansatz} \label{sec:imprAnsatz}

We construct the improved approximation in system~\eqref{modeq4aneu}, which is of the form
\begin{align}
\begin{split}
\partial_t Z_1 & = \lambda_+(-\i\partial_x) Z_1 + \theta(\partial_x) g_1(Z_1,Z_2),\\
\partial_t Z_2 & = \lambda_-(-\i\partial_x) Z_2 +  g_2(Z_1,Z_2),
\end{split} \label{modeq4form}
\end{align}
with $ g_1 $, $ g_2 $ being smooth   nonlinear terms
for which we can compute their Taylor expansion up to arbitrary order (see   Appendix~\ref{appA2}).
We  choose $ c $  as in~\eqref{veloc2} and
make the ansatz
\begin{align}
\begin{split}
 Z_1(x,t) & =  \varepsilon^2 A_0(\varepsilon x , \varepsilon^3 t ) + \ldots + \varepsilon^{N+2} A_N(\varepsilon x , \varepsilon^3 t )  , \\
  Z_2(x,t) & =  \varepsilon^4 B_0(\varepsilon x , \varepsilon^3 t ) + \ldots + \varepsilon^{N+3} B_{N-1}(\varepsilon x , \varepsilon^3 t ) .
 \end{split} \label{imprAnsatz3}
\end{align}
Under the scaling $\xi = \varepsilon x$,
using the expansions of $\lambda_{\pm}(k)$ about $k = 0$ from $~\S\ref{sec2} $,
we have the following formal expansion
in the parameter regime~\eqref{parreg}
\begin{align*}
 \lambda_+(-\i\partial_x) & = \varepsilon^3 \gamma_{lin} \partial_{\xi}^3 + \mathcal{O}(\varepsilon^4),\\
  \lambda_-(-\i\partial_x) & = - 2 \sigma_s + \mathcal{O}(\varepsilon),\\
  \theta(\partial_x) & = \varepsilon \partial_{\xi}
\end{align*}
with $ \gamma_{lin} = -c_{3,s} $ as in~\S\ref{sec3}. By substitution of~\eqref{imprAnsatz3} into~\eqref{modeq4form}, we therefore find at $ \mathcal{O}(\varepsilon^5)$ in the first component of~\eqref{modeq4form} the KdV equation
\begin{equation} \label{kdvcec}
\partial_\tau A_0 = \gamma_{lin} \partial_{\xi}^3 A_0  + \gamma_{non} \partial_{\xi}(A_0^2),
\end{equation}
with $\gamma_{non} = \beta - \alpha$ as in~\S\ref{sec3}. At $ \mathcal{O}(\varepsilon^4) $ in the second component of~\eqref{modeq4form} we find a linear
algebraic equation
\begin{align*}
0 & = -2 \sigma_s B_0 + c_* (A_0)^2 ,
\end{align*}
where here and in the following  various coefficients are denoted by the  symbol $ c_* $.
\begin{remark}{\label{compareAnsatz} \rm
We computed for the ansatz in~\S\ref{sec3}  that $B = - \sigma^{-1} A + \ord(\epsilon)$, such that
$$
V_{app}^{*,\epsilon} = \epsilon^2 \left(\begin{array}{c} A \\ B \end{array} \right) + \mathcal{O}(\epsilon^3)
=  \epsilon^2  A \left(\begin{array}{c} 1 \\ - \sigma^{-1}   \end{array} \right) + \mathcal{O}(\epsilon^3).
$$
On the other hand we have
\begin{align*}
V_{app}^{\kappa,\epsilon} &= S_{\vartheta}^{-1} S_{diag}^{-1} Z = \epsilon^2\left(\begin{array}{cc} 1 & * \\ - \sigma^{-1} & * \end{array} \right) A_0 \left(\begin{array}{c} 1 \\ 0  \end{array} \right) + \mathcal{O}(\epsilon^3)\\ &= \epsilon^2  A_0 \left(\begin{array}{c} 1 \\ - \sigma^{-1}   \end{array} \right) + \mathcal{O}(\epsilon^3).
\end{align*}
{Therefore, we  have formally $A = A_0 + \mathcal{O}(\varepsilon)$,  and so
the coefficients in~\eqref{eq:kdv} and~\eqref{kdvcec} are (to leading order) the same. On the other hand, we note that $B$ and $B_0$ are not equal.}
}\end{remark}
At
$ \mathcal{O}(\varepsilon^6) $ in the first component of~\eqref{modeq4form} we find a linearized inhomogeneous
KdV equation
\begin{align*}
\partial_{\tau} A_1 & =   \gamma_{lin}  \partial_{\xi}^3 A_1  +  2  \gamma_{non}  \partial_{\xi} (A_0 A_1) + f_{A,1}(A_0) ,
\end{align*}
where $  f_{A,1}(A_0)  $ is a function, solely depending on $ A_0 $ and its $\xi$-derivatives, given by
$$
 f_{A,1}(A_0)  = c_* \partial_{\xi}^2 A_0 + c_* \partial_{\xi}^4 A_0 + c_* \partial_{\xi}^2 ((A_0)^2).
$$
The principal structure of the subsequent equations remains the same. For example
at $ \mathcal{O}(\varepsilon^5) $ in the second component of~\eqref{modeq4form} we find a linear
algebraic equation
\begin{align*}
0 & = -2 \sigma_s B_1 +  f_{B,1}(A_0,A_1,B_0),
\end{align*}
and at
$ \mathcal{O}(\varepsilon^7) $ in the first component of~\eqref{modeq4form} we find a linearized inhomogeneous
KdV equation
\begin{align*}
\partial_{\tau} A_2 & = \gamma_{lin} \partial_{\xi}^3 A_2  +  2 \gamma_{non} \partial_{\xi} (A_0 A_2) + f_{A,2}(A_0,A_1,B_0) ,
\end{align*}
where $  f_{B,1} $ and $f_{A,2}$ are functions which  solely depend on
the solutions to the equations before, i.e., here on $ A_0,A_1 $, and $ B_0$ and their $\xi$- and $\tau$-derivatives.
For $ m\in \{3,\ldots,N\} $ at $ \mathcal{O}(\varepsilon^{m+3}) $ in the second component of~\eqref{modeq4form} we find a linear
algebraic equation
\begin{align*}
0 & = -2 \sigma_s B_{m-1} +  f_{B,m-1}(A_0,\ldots,A_{m-1},B_0,\ldots,B_{m-2}),
\end{align*}
and at
$ \mathcal{O}(\varepsilon^{m+5}) $ in the first component of~\eqref{modeq4form} we find a linearized inhomogeneous
KdV equation
\begin{align*}
\partial_{\tau} A_m & = \gamma_{lin} \partial_{\xi}^3 A_m  +  2 \gamma_{non} \partial_{\xi} (A_0 A_m) + f_{A,m}(A_0,\ldots,A_{m-1},B_0,\ldots,B_{m-2}) ,
\end{align*}
where again $  f_{B,m-1} $ and $f_{A,m} $ are functions which  solely depend on
the solutions $A_0,\ldots,A_{m-1}, B_0,\ldots,B_{m-2}$ and their $\xi$- and $\tau$-derivative.
In this procedure also the temporal derivatives of $B_m$-terms occur. They can be expressed
in terms of spatial derivatives of the solutions to the equations before by differentiating
the algebraic $ B_m $-equation with respect to time and then expressing the temporal derivatives of the solutions to the equations before by the right hand sides.
For instance we find
$$
0  = -2 \sigma_s  \partial_{\tau} B_0 + c_* A_0  \partial_{\tau} A_0 =  - 2 \sigma_s  \partial_{\tau} B_0+ c_* A_0
(\gamma_{lin} \partial_{\xi}^3 A_0 +  \gamma_{non} \partial_{\xi} ((A_0)^2)).
$$
Hence, if local existence and uniqueness as well as sufficient regularity can be established for these equations,
we can solve them step by step.
Before we go on, we remark that $  f_{B,m-1} $ and $f_{A,m} $ only contain finitely many derivatives.

\subsection{Local existence and uniqueness for the extended amplitude system}
\label{sec:imprAnsatz2}

In this section we establish local existence and uniqueness  in  $H^\infty_{\mu,s}$-spaces for the
extended amplitude system

\begin{align}
\begin{split}
\partial_{\tau} A_0  & =  \gamma_{lin} \partial_{\xi}^3 A_0 +  \gamma_{non} \partial_{\xi} ((A_0)^2), \\
\partial_{\tau} A_1 & = \gamma_{lin} \partial_{\xi}^3 A_1  +  2 \gamma_{non} \partial_{\xi} (A_0 A_1) + f_{A,1}(A_0) , \\
   & \vdotswithin{=}\\
\partial_{\tau} A_m & = \gamma_{lin} \partial_{\xi}^3 A_m  +  2 \gamma_{non} \partial_{\xi} (A_0 A_m) + f_{A,m}(A_0,\ldots,A_{m-1},B_0,\ldots,B_{m-2}) ,
\\
   & \vdotswithin{=} \\
\partial_{\tau} A_N & = \gamma_{lin} \partial_{\xi}^3 A_N  +  2 \gamma_{non} \partial_{\xi} (A_0 A_N) + f_{A,N}(A_0,\ldots,A_{N-1},B_0,\ldots,B_{N-2}),
\end{split}\notag
\intertext{with}
\begin{split}
0 & = - 2 \sigma_s  B_0 + c_* (A_0)^2 ,  \\
  & \vdotswithin{=}  \\
0 & = - 2 \sigma_s B_{m-1} +  f_{B,m-1}(A_0,\ldots,A_{m-1},B_0,\ldots,B_{m-2}),\\
  & \vdotswithin{=} \\
0 & = - 2 \sigma_s B_{N-1} +  f_{B,N-1}(A_0,\ldots,A_{N-1},B_0,\ldots,B_{N-2}),
\end{split}\label{algebraiceq}
\end{align}
which we derived in~\S\ref{sec:imprAnsatz}.

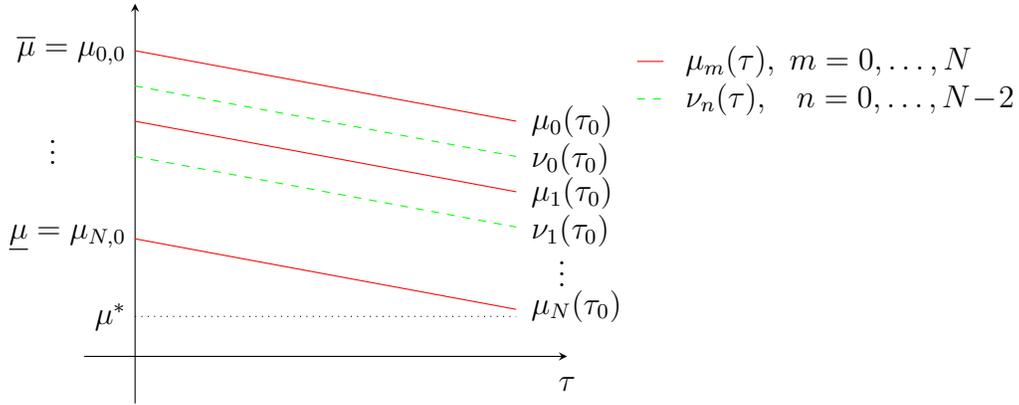
\begin{figure}[htbp] 
   \centering
  \begin{tikzpicture}[rotate=0,scale=1]
    \begin{axis}[
	xmin=-0.2, xmax=1.7,
	ymin=-0.2, ymax=1.5,
	axis lines=center,
	ticks=none,
	xlabel={$ \tau $},
	x label style={at={(axis description cs:1.0,0.05)},anchor=center},
	width=8.cm,
    ]
      \addplot+[red, domain=0:1.5,samples=5,no marks] {1.3-0.2*x};
      \addplot+[green,dashed, domain=0:1.5,samples=5,no marks] {1.15-0.2*x};
      \addplot+[red, domain=0:1.5,samples=5,no marks] {1.0-0.2*x};
      \addplot+[green,dashed, domain=0:1.5,samples=5,no marks] {0.85-0.2*x};
      \addplot+[red, domain=0:1.5,samples=5,no marks] {0.5-0.2*x};
      \addplot+[black,dotted, domain=0:1.5,samples=5,no marks] {0.17};
    \end{axis}
    \node[black,anchor=east] at (0.7,4.7) {$\overline \mu = \mu_{0,0}$};
    \node[black,anchor=east] at (-0.215,3.45) {$\vdots$};
    \node[black,anchor=east] at (0.7,2.2) {$\underline \mu = \mu_{N,0}$};
    \node[black,anchor=east] at (0.7,1.15) {$\mu^*$};
    \node[black,anchor=west] at (5.8,3.75) {$\mu_{0}(\tau_0)$};
    \node[black,anchor=west] at (5.8,3.25) {$\nu_{0}(\tau_0)$};
    \node[black,anchor=west] at (5.8,2.8) {$\mu_{1}(\tau_0)$};
    \node[black,anchor=west] at (5.8,2.30) {$\nu_{1}(\tau_0)$};
    \node[black,anchor=west] at (6.15,1.83) {$\vdots$};
    \node[black,anchor=west] at (5.8,1.3) {$\mu_{N}(\tau_0)$};
    \draw[red] (7.35,4.55) -- (7.70,4.55);
    \draw[green,dashed] (7.35,4.05) -- (7.70,4.05);
    \node[black,anchor=west] at (7.85,4.55) {$\mu_{m}(\tau), \; m=0,\dots, N$};
    \node[black,anchor=west] at (7.85,4.05) {$\nu_{n}(\tau), \!\quad n=0,\dots, N\!-\!2$};
  \end{tikzpicture}
   \caption{The left panel shows a lower bound to the width of the strip of analyticity  of the functions $ A_0, B_0, A_1, \dotsc, B_{N-2}, A_N $
  as a function of time $ \tau $.
}
   \label{figregularity0}
\end{figure}

In order to do so we take three numbers $ \overline{\mu}, \underline{\mu} $ and $\mu_*$ with  $ \overline{\mu} > \underline{\mu} > \mu_{*}>0$.
Let
 $$ \mu_{m,0} = \overline{\mu}  - (\overline{\mu} -\underline{\mu}) \frac{m}{N}
 \qquad \textrm{and} \qquad
 \mu_m(\tau) =  \mu_{m,0} - \eta \tau
 $$
 for $0 < \eta \tau < \underline{\mu} - \mu_{*} $, where $\eta > 0$ is some $\epsilon$-independent constant.
 Then we choose $ A_m|_{\tau=0} \in H^\infty_{\mu_{m,s}}$, $s$ sufficiently large, and look for solutions with
  $ A_m(\tau) \in H^\infty_{\mu_{m}(\tau),s}$.
  Moreover, we set $ \nu_m(\tau) = \mu_m(\tau) - \frac{1}{2N}  (\overline{\mu} -\underline{\mu}) $
and look for solutions with
  $ B_m(\tau) \in H^\infty_{\nu_{m}(\tau),s}$ -- see Figure~\ref{figregularity0}.

Since $ f_{A,m}=  f_{A,m}(A_0,\ldots,A_{m-1},B_0,\ldots,B_{m-2}) $
only contains finitely many derivatives, we
have for $s \geq 0$ that
\begin{equation*}
f_{A,m}: H^\infty_{\mu_{0}(\tau),s} \times \ldots \times H^\infty_{\mu_{m-1}(\tau),s}
\times H^\infty_{\nu_{0}(\tau),s} \times \ldots \times H^\infty_{\nu_{m-2}(\tau),s}
\to H^\infty_{\mu_{m}(\tau),s}
\end{equation*}
is a bounded mapping.
Similarly, we have that
\begin{equation*}
f_{B,m-1}: H^\infty_{\mu_{0}(\tau),s} \times \ldots \times H^\infty_{\mu_{m-1}(\tau),s}
\times H^\infty_{\nu_{0}(\tau),s} \times \ldots \times H^\infty_{\nu_{m-2}(\tau),s}
\to H^\infty_{\nu_{m-1}(\tau),s}
\end{equation*}
is a bounded mapping. Hence, substituting $ B_m $ in $ f_{A,m} $ by
$ f_{B,m}/(2 \sigma_s) $
gives a system of the form
\begin{align} \label{cfb1}
\partial_{\tau} A_0  & =  \gamma_{lin} \partial_{\xi}^3 A_0 +   \gamma_{non}  \partial_{\xi} (A_0^2), \\
 \label{cfb2}
\partial_{\tau} R  & = \gamma_{lin} \partial_{\xi}^3 R +  2 \gamma_{non}  \partial_{\xi} (A_0 R) + F(A_0 ,R),
\end{align}
with  $ R=(A_1,\ldots,A_{N}) $, where $ F:H^\infty_{\mu_{0}(\tau),s} \times \h^\infty_{\tau,s} \to \h^\infty_{\tau,s} $, with
$$
\h^\infty_{\tau,s} = H^\infty_{\mu_{1}(\tau),s} \times \ldots \times H^\infty_{\mu_{N}(\tau),s}
$$
is a smooth mapping.

We
introduce $ \widehat{\rho}_m(k,\tau) = e^{\mu_m(\tau)|k|} $ and $ {\A}_m(\xi,\tau) $
via
\begin{align}
\widehat{\A}_m(K,\tau) = \widehat{\rho}_m(K,\tau)  \widehat{A}_m(K,\tau). \label{anasmooth}
\end{align}
Moreover,  we introduce $ {\cR} = ({\A}_1,\ldots,{\A}_{N})$ and
$ \varrho_N(\tau) = \textrm{diag}({\rho}_1(\tau),\ldots,{\rho}_{N}(\tau)) $ and
have the following  local existence and uniqueness results
for the KdV equation and for the extended amplitude system.

\begin{theorem} \label{linkin1}
Let $ A_0|_{\tau = 0} \in  H^\infty_{\mu_{0}(0),s} $ and $s\geq 3$. Then there exists $\tau_0 > 0 $ and a solution
$ A_0 $ to~\eqref{cfb1} for $ \tau \in [0,\tau_0] $ with  $ A_0(\tau) \in H^\infty_{\mu_{0}(\tau),s}$,
 where $0 < \eta \tau_0 < \underline{\mu} - \mu_{*} $. Moreover
 $\A_0 \in C([0,\tau_0],H^{s}) \cap C^1([0,\tau_0],L^2(\R))$.
\end{theorem}
\begin{theorem} \label{th66}
Let $ ( A_0|_{\tau = 0} , R|_{\tau = 0} ) \in  H^\infty_{\mu_{0}(0),s} \times  \h^\infty_{0,s}$ and $s\geq 3$. Then there exists $ \tau_0 > 0 $ and solutions
$ (A_0,R) $ to~\eqref{cfb1}-\eqref{cfb2} for $ \tau \in [0,\tau_0] $
with  $ (A_0,R)(\tau) \in H^\infty_{\mu_{0}(\tau),s} \times \h^\infty_{\tau,s}$,
 where $ 0 < \eta \tau_0 < \underline{\mu} - \mu_{*} $.
 Moreover, $(\A_0,\cR) \in C([0,\tau_0],H^{s}) \cap C^1([0,\tau_0],L^2(\R))$.
 \end{theorem}
\noindent
{\bf Proof.}
In order to prove local existence and uniqueness
of solutions to~\eqref{cfb1} and~\eqref{cfb1}-\eqref{cfb2}
and to
estimate  the $ H^\infty_{\mu_{0}(\tau)} $-norm and the $ H^\infty_{\mu_{0}(\tau)} \times \h^\infty_{\tau} $-norm of these solutions, we proceed as in~\S\ref{sec:diagonal} by introducing time-dependent exponential weights via~\eqref{anasmooth}. This transforms~\eqref{cfb1}-\eqref{cfb2} into the system
\begin{align} 
\begin{split}
\partial_{\tau} {\A}_0(\tau)  & = - \eta |K|_{op}  {\A}_0(\tau)+  \gamma_{lin} \partial_{\xi}^3 {\A}_0(\tau) \\
			      & \quad +   \gamma_{non} \rho_0(\tau) \partial_{\xi} ( ( \rho_0^{-1}(\tau){\A}_0(\tau))
( \rho_0^{-1}(\tau)  {\A}_0(\tau))) , \end{split}\label{cfb1s}
\\
\begin{split}
\partial_{\tau} {\cR}(\tau)  & = - \eta |K|_{op} {\cR}(\tau)  +
 \gamma_{lin} \partial_{\xi}^3 {\cR}(\tau) \\
  & \quad +  2 \gamma_{non}  \varrho_N(\tau) \partial_{\xi} (( \rho_0^{-1}(\tau){\A}_0(\tau))( \varrho_N^{-1}(\tau)  {\cR}(\tau)))
  \\
  & \quad + \check{F}({\A}_0(\tau)) ,{\cR}(\tau)),
  \end{split} \label{cfb2s}
\end{align}
with $ \widehat{  |K|_{op}  {\A}_0} (K)=  |K| \widehat{{\A}}_0(K) $ and
where $ \check{F}({\A}_0 ,{\cR}) $ is the transformed nonlinearity. Note that $\check{F}$ has a strict upper triangular structure, i.e. $\check F_i$ only depends on $\A_1, \dots, {\A}_{i-1}$. Moreover, if ${\A}_0, \dots, {\A}_{i-1} \in C([0,\tau_0],H^{s})$ then $\check{F}_{i} \in C([0,\tau_0],H^{s})$, $1 \leq i \leq N$.

System~\eqref{cfb1s}-\eqref{cfb2s} is a  classical
quasilinear system in sense of~\cite{Kato75}. For~\eqref{cfb1s} and~\eqref{cfb1s}-\eqref{cfb2s} we obtain local existence and uniqueness
in the space $C([0,\tau_0],H^{s}) \cap C^1([0,\tau_0],L^2(\R))$ for $ s \geq 3 $.
We remark that since the quasilinear parts of~\eqref{cfb1s} and~\eqref{cfb1s}-\eqref{cfb2s}
have the same structure the time of existence of the solutions, which we obtain from our construction, is the same  for~\eqref{cfb1s}  and~\eqref{cfb1s}-\eqref{cfb2s}.
 \qed

\begin{remark}{\rm
Alternatively, one can establish local existence and uniqueness via optimal regularity results. Since the linear semigroup gains one derivative and since the nonlinear terms lose
one derivative, optimal regularity results in the sense of~\cite{LM1,LM2}
can be used. We remark that the solutions to the KdV equation stay analytic for all times~\cite{katoMasuda1986,tarama2004,bona}.
}\end{remark}

%
%
%
%
%
%
%
%


Starting from the solution $A \in C([0,\tau_0],H^\infty_{\mu_{\A}})$ to the KdV equation~\eqref{KdV} in Theorem~\ref{mainresult7} and taking $\mu_* < \mu_0(0) < \mu_\A$ and $A_0|_{\tau = 0} = A|_{\tau=0}$ in Theorem~\ref{th66}, we obtain $\tilde{\tau}_0 \in (0,\tau_0]$ and an improved approximation $(A_0,A_1,\ldots,A_{\lceil\kappa\rceil})(\tau)$ lying in $H^\infty_{\mu_*,s}$ such that $A(\tau) = A_0(\tau)$ for $\tau \in [0,\tilde{\tau}_0]$. Now, define $Z_{app}^{\kappa,\epsilon}$ by~\eqref{imprAnsatz3}, where the $B_i$ are defined through the algebraic equations~\eqref{algebraiceq}. As $Z_{app}^{\kappa,\epsilon}(t)$ lies in $H^{\infty}_{\mu_*/\epsilon,s}$ for $0 \leq \epsilon^3 t \leq \tilde \tau_0$, the transformed improved approximation $\cZ_{app}^{\kappa,\epsilon}(t) := S_\omega(t) Z_{app}^{\kappa,\epsilon}(t)$ satisfies the desired estimate~\eqref{impresolvz} in Theorem~\ref{mainresult7} by construction (taking $s > 0$ larger if necessary), where we use $\mu(t) \geq \mu_*/\epsilon$ for $0 \leq \epsilon^3 t \leq \tau_0$. All that remains to show is the error estimate~\eqref{Hmboundz}, which will be established in~\S\ref{sec:error}.

\subsection{The residual estimates} \label{sec:resest}

For the improved approximation $ Z_{app}^{\kappa,\epsilon} $ defined in~\S\ref{sec:imprAnsatz}-\ref{sec:imprAnsatz2} we have  the following estimate on the residual terms.
\begin{lemma}
There exists $C_{res} > 0$ such that for all $\epsilon \in (0,1)$ we have
\begin{align}
\sup_{0 \leq \epsilon^{3} t \leq \tilde{\tau}_0}\|\Res_{Z}(Z_{app}^{\kappa,\epsilon}(t))\|_{H^\infty_{\mu_*/\epsilon}} \leq C_{res}\epsilon^{3+\kappa},
\label{resbB}
\end{align}
where we fixed $ \mu_* \in (0,\mu_{\A}) $ in~\S\ref{sec:diagonal} and $\Res_Z$ is defined in Remark~\ref{remresz}.
\end{lemma}
\noindent
{\bf Proof.}
By the construction of the improved approximation
in~\S\ref{sec:imprAnsatz} we eliminated formally all terms of order
$ \mathcal{O}(\varepsilon^{3+\kappa+1/2}) $ in the residual defined in~\eqref{resform}.
In~\S\ref{sec:imprAnsatz2} we showed that for any fixed $s > 0$
the improved approximation $Z_{app}^{\kappa,\epsilon}(t)$ is in $ H^{\infty}_{\mu_*/\epsilon,s} $ and so are the remaining terms
small in $ H^{\infty}_{\mu_{*}/\epsilon,m} $ (after taking $s$ sufficiently large).
In detail we use Lemma~\ref{applem4} and apply it on $ \lambda_+(k) + i \gamma_{lin} k^3 = \mathcal{O}(|k|^4) $ as $k \to 0$. A straightforward multilinear generalization of
Lemma~\ref{applem4}  is applied  on the associated  kernels appearing in the representation of the
nonlinear terms in Fourier space. As an example, in Fourier space bilinear terms can be written
as
$$
\int b(k,k-m,m) \widehat{u}(k-m)\widehat{u}(m) dm
$$
with kernel $ b $, which is smooth at the origin and can be  expanded like $ \lambda_{\pm} $.
\qed

\subsection{The error estimates} \label{sec:error}

Our starting system  is equation~\eqref{modeq4c}
\begin{align}  \label{modeq4d}
\partial_t \cZ = L_\cZ \cZ + N_\cZ(t,\cZ),
\end{align}
with linearity $L_\cZ$ having the Fourier symbol
$ \widehat{L_{\cZ}}(k) = \mathrm{diag}\left(\lambda_{2,+}(k),\lambda_{2,-}(k)\right) $,
and where the nonlinear terms are by Lemma~\ref{applem3} of the form
\begin{align*}
N_\cZ(t,\cZ) & =
 \begin{pmatrix}  \theta  (B_1(t,\cZ_1,\cZ_1)+ B_2(t,\cZ_1,\cZ_2)+ B_3(t,\cZ_2,\cZ_2) +  \h_1(t,\cZ)) \\ B_4(t,\cZ_1,\cZ_1)+ B_5(t,\cZ_1,\cZ_2)+ B_6(t,\cZ_2,\cZ_2) +  \h_2(t,\cZ) \end{pmatrix} ,
\end{align*}
in which the bilinear terms $B_j$ enjoy the estimate
\begin{align}
\|B_j(t,\cZ,\cW)\|_{H^s} \leq C\|\cZ\|_{H^{s+1}}\|\cW\|_{H^{s+1}}, \qquad j = 1,\ldots,6, \label{estbil}
\end{align}
and where the higher order terms $ \h_{1,2} $ obey
\begin{align*}
\| \h_{1,2}(t,\cZ) \|_{H^s} \leq C \|\cZ\|_{H^{s}}^2\|\cZ\|_{H^{s+1}},
\end{align*}
for small $ \|\cZ\|_{H^{s+1}} $.
\begin{remark}
{\rm
System~\eqref{modeq4d} has the properties of a  semilinear parabolic system and so local existence and uniqueness is clear~\cite{Lunardi}. The solutions exist as long as we can bound them with our subsequent error estimates.
}
\end{remark}
 The improved long-wave KdV approximation
$Z_{app}^{\kappa,\epsilon}(t)$ defined in~\S\ref{sec:imprAnsatz2} for the variable
$ Z $ gives via
\begin{align}
\cZ_{app}^{\kappa,\epsilon}(t) = S_{\omega}(t) Z_{app}^{\kappa,\epsilon}(t) =: \varepsilon^2 \cW(t) = \begin{pmatrix}  \varepsilon^2 \cW_1(t) \\ \varepsilon^4 \cW_2(t) \end{pmatrix} \label{strucZapp}
\end{align}
an improved long-wave KdV approximation for the variable $ \cZ $.

We introduce the error function $ R $ by
\begin{equation*}
\cZ = \varepsilon^2 \cW + \varepsilon^{\kappa} R,
\end{equation*}
with $R=(R_1,R_2)$.
The error function $ R $ satisfies
\begin{align*}
 \partial_t R = L_\cZ R + \Non(t,R) +  \epsilon^{-\kappa}  \Res_{\cZ}(\varepsilon^2 \cW),
\end{align*}
with $ R|_{t=0} = 0 $, where
\begin{align}
\begin{split}
{\Non}(t)(R) &:= \epsilon^{-\kappa}  [{N_\cZ}\left(t,\varepsilon^2 \cW+ \epsilon^\kappa R) - {N_\cZ}(t,\varepsilon^2 \cW)\right] =: \begin{pmatrix} {\theta}{\Non}_1(t)(R) \\ {\Non}_2(t)(R)\end{pmatrix}.
\end{split} \label{strucG2}
 \end{align}
We have the  expansion
\begin{equation} \label{expansion}
 \begin{pmatrix} {\theta} {\Non}_1(t)(R) \\ {\Non}_2(t)(R)\end{pmatrix}
 =  \begin{pmatrix} \varepsilon^2 {\theta} (2 B_1(t,\cW_1,R_1) + B_2(t,\cW_1,R_2))
 \\  \varepsilon^2  (2 B_4(t,\cW_1,R_1) + B_5(t,\cW_1,R_2) 
 )  \end{pmatrix}
 +  \begin{pmatrix} \h_3(t,R)
 \\ \h_4(t,R)
  \end{pmatrix},
\end{equation}
where $\h_{3,4}$ obey by Corollaries~\ref{appcoro1} and~\ref{applem5} and Lemma~\ref{applem3} the estimate
\begin{align}
\| \h_{3,4}(t,R) \|_{H^m} \leq C  \varepsilon^4 \| R \|_{H^{m+1}}
+ C_2(M) \varepsilon^{\kappa} \| R \|_{H^{m}} \| R \|_{H^{m+1}} , \label{esthih2}
\end{align}
as long as $ \| R \|_{H^{m}}  \leq M $, with $ M $ defined below. More precisely, one observes from~\eqref{strucG2} and~\eqref{expansion} that the terms in $\h_{3,4}(t,R)$ that are nonlinear in $R$ are estimated by $C_2(M) \varepsilon^{\kappa} \| R \|_{H^{m}} \| R \|_{H^{m+1}}$,
where $ C_2(M) $ denotes a constant depending on $ M $. On the other hand, the terms in $\h_{3,4}(t,R)$ that are linear in $R$ are either bilinear terms of the form $B_j(t,\cW_2,R_i)$ or $B_j(t,R_i,\cW_2)$ with $j = 1,\ldots,6$, $i = 1,2$ or they are (at least) quadratic in $\cW$. Consequently, these terms are estimated by $C \varepsilon^4 \| R \|_{H^{m+1}}$ using~\eqref{strucZapp}.

The error is estimated by the energy
$$
\E(t) := \frac{1}{2}(\|R_1(t)\|_{H^m}^2+\|R_2(t)\|_{H^m}^2) = \frac{1}{2} \|R(t)\|_{H^m}^2 .
$$
We compute
\begin{align}
 \partial_t \E = \Re\left(\langle R, L_\cZ R \rangle_{H^m}+ \langle R,\Non(R) \rangle_{H^m}+
 \langle R,  \epsilon^{-\kappa}  \Res_Z(\varepsilon^2 W) \rangle_{H^m}\right)\label{erroreq22}.
\end{align}

{\bf a)} We start with the bound on  the linear term $\langle R,L_\cZ R \rangle_{H^m}$ in~\eqref{erroreq22}. Using that $L_\cZ$ has the Fourier symbol~\eqref{Fsymb} satisfying the bounds~\eqref{eigvbounds}, we obtain the existence of
a constant $c_1 > 0$ such that
\begin{align*}  
  \Re\left(\langle R,L_\cZ R\rangle_{H^m}\right) \leq - \epsilon^2 \frac{\eta}{2} \||k|_{op}^{1/2} R_1\|_{H^{m}}^2 - c_1\|R_2\|_{H^{m+1}}^2,
\end{align*}
is satisfied for  $\varepsilon > 0$ sufficiently small, where, as above, $ |k|_{op}^{1/2} R_j $
is defined in Fourier space by $\fourier (|k|_{op}^{1/2} R_j)(k) = |k|^{1/2} \widehat{R}_j(k)$.

{\bf b)} Next, we bound the residual term $\langle R, \Res_{\cZ}(\varepsilon^2 \cW)\rangle_{H^m}$ in~\eqref{erroreq22}  by employing the estimate~\eqref{resbB}. Thus, using that $\mu(t) \leq \mu_*/\epsilon$ and $S_{\omega}(t)$ maps  $H^\infty_{\mu(t)}$ continuously into $H^m$, it holds
\begin{align*}
  \begin{split}
  \left|\langle R, \Res_\cZ(\varepsilon^2 \cW)\rangle_{H^m}\right| &\leq C\|R\|_{H^m}\|\Res_Z\|_{H^\infty_{\mu(t)/\epsilon}} \\&\leq C\epsilon^{3+\kappa}\left(1+\E\right),
  \end{split}
\end{align*}
where we used the Cauchy-Schwarz and  Young's inequalities.

{\bf c)} Finally, we bound the nonlinear term
\begin{equation*}
  \left\langle R, \Non(t,R)\right\rangle_{H^m} = \left\langle R_1, \theta \Non_1(t,R)\right\rangle_{H^m}  +  \left\langle R_2,  \Non_2(t,R)\right\rangle_{H^m} .
\end{equation*}

{\bf c1)}
Using Corollaries~\ref{appcoro1} and~\ref{applem5},~\eqref{estbil},~\eqref{expansion} and~\eqref{esthih2}
gives
\begin{align*}
\begin{split}
&\left|\left\langle R_2, \Non_2(t,R) \right\rangle_{H^m}\right|
\\
& \qquad \leq \left|\left\langle R_2, \Non_2(t,R) \right\rangle_{H^{m-1}} \right| + \left|\left\langle |k|_{op}^{1/2}R_2, |k|_{op}^{1/2} \Non_2(t,R) \right\rangle_{H^{m-1}}\right|
\\
& \qquad \leq 2 \|R_2\|_{H^{m+1/2}} \|{\Non}_2(t,R)\|_{H^{m-1/2}}\\
& \qquad \leq \|R_2\|_{H^{m+1/2}} \left(C \epsilon^2  \|R_2\|_{H^{m+1/2}} + C \epsilon^2 \|R_1\|_{H^{m+1/2}} \right. \\ & \qquad \qquad  \qquad \qquad \qquad \qquad \left. + C_2(M)  \epsilon^{\kappa}  \|R \|_{H^{m}}  \|R \|_{H^{m+1/2}} \right)
\\
& \qquad \leq C \epsilon \|R_2\|_{H^{m+1/2}}^2 + C \epsilon^3 \|R\|_{H^{m+1/2}}^2
+ C_2(M)  \epsilon^{\kappa}  \|R \|_{H^{m}}  \|R \|_{H^{m+1/2}}^2
 \end{split}
\end{align*}
under the assumption $ \|R\|_{H^{m}} \leq M $,
where we used $ \epsilon^2 a b \leq  \epsilon a^2 +  \epsilon^3  b^2 $ in the last line.

{\bf c2)}
Similarly, we estimate
\begin{align*}
\begin{split}
&\left|\left\langle R_1, \h_3(t,R) \right\rangle_{H^m}\right| \leq 2 \|R_1\|_{H^{m+1/2}} \|{\h}_3(t,R)\|_{H^{m-1/2}}\\
& \qquad \leq \|R_1\|_{H^{m+1/2}}  \left(C \epsilon^4 \|R\|_{H^{m+1/2}} + C_2(M) \epsilon^{\kappa}  \|R \|_{H^{m}}  \|R \|_{H^{m+1/2}}\right )
\\
& \qquad \leq  \left(C \epsilon^4 + C_2(M)\epsilon^{\kappa}\right) \|R \|_{H^{m+1/2}}^2 ,
 \end{split}
\end{align*}
as long as $ \|R\|_{H^{m}} \leq M $.

\newpage
{\bf c3)}
Finally, using that the symbol of $ |k|_{op}^{-1} \theta $ is bounded at the wave number $ k = 0 $, we find by~\eqref{estbil}
\begin{align*}
\begin{split}
&\left|\left\langle R_1, \varepsilon^2  \theta B_j(t,\cW_1,R_j) \right\rangle_{H^m}\right| \leq \epsilon^2 \||k|_{op}^{1/2} R_1\|_{H^m} \||k|_{op}^{-1/2}  \theta B_j(t,\cW_1,R_j) \|_{H^m}\\
\qquad & \leq  \varepsilon^2  \||k|_{op}^{1/2}R_1\|_{H^m} \|\theta^{1/2}  B_j(t,\cW_1,R_j) \|_{H^{m-1/2}}\\
\qquad & \leq C \epsilon^2 \||k|_{op}^{1/2}R_1\|_{H^m} \|\theta^{1/2}   \cW_1 \|_{W_{m+1/2}} \|R\|_{H^{m+1/2}} \\
\qquad & \quad + C \epsilon^2 \||k|_{op}^{1/2}R_1\|_{H^m} \| \cW_1 \|_{W_{m+1/2}} \|\theta^{1/2} R \|_{H^{m+1/2}} \\
\qquad & \leq C  \epsilon^{5/2}  \||k|_{op}^{1/2}R_1\|_{H^m}  \|R\|_{H^m} +C \epsilon^2 \||k|_{op}^{1/2}R_1\|_{H^m} \||k|_{op}^{1/2} R \|_{H^m}\\
\qquad & \leq C (\epsilon^2 \||k|_{op}^{1/2}R_1\|_{H^m}^2 +  \epsilon^{2}  \|R_2\|_{H^{m+1/2}}^2+\epsilon^3  \|R_1 \|_{H^{m+1/2}}^2)
\end{split}
\end{align*}
where we used that $ \|\theta^{1/2}   \cW_1 \|_{W_{m+1/2}} = \mathcal{O}(\varepsilon^{1/2}) $ due to Corollary~\ref{appcoro4}, and $ \epsilon^{5/2} a b \leq  \epsilon^2 a^2 +  \epsilon^3  b^2 $. The bound on $\left|\left\langle R_1, \varepsilon^2  \theta B_j(R_j,\cW_1) \right\rangle_{H^m}\right|$ is analogous.

{\bf d)} Before we summarize all estimates
we introduce the following notation. Constants which do not depend on $ \E $ are denoted
by $ C_1 $, constants only  depending on the residual are called $ C_{res} $, and constants
which depend on $ M $ are called, with some slight abuse of notation, again $ C_2(M) $.

Applying the above estimates to~\eqref{erroreq22} yields
\begin{align*}
\partial_t \E & \leq - \epsilon^2 \frac{\eta}{2} \|K^{\tfrac{1}{2}} R_1\|_{H^{m}}^2 - c_1\|R_2\|_{H^{m+1}}^2
 \\ &\quad  + C_1 (\epsilon^2 \||k|_{op}^{1/2}R_1\|_{H^m}^2
+  \epsilon^{2}  \|R_2\|_{H^{m+1/2}}^2+\epsilon^3  \|R_1 \|_{H^{m+1/2}}^2)
 \\ &\quad  +
C_1 \epsilon \|R_2\|_{H^{m+1/2}}^2 + C_1 \epsilon^3 \|R\|_{H^{m+1/2}}^2
+ C_2(M) \epsilon^{\kappa}  \|R \|_{H^{m}}  \|R \|_{H^{m+1/2}}^2
 \\ &\quad
+C_{res}\epsilon^{3}\left(1+\E(t)\right)
\\ & \leq (- \epsilon^2 \frac{\eta}{2} +  C_1 \epsilon^2 +  C_1 \epsilon^3 + C_2(M) \epsilon^{\kappa}  \|R \|_{H^{m}} ) \||k|_{op}^{1/2} R_1\|_{H^{m}}^2
\\ &\quad
+ (- c_1+ C_1 \epsilon^2 + C_1 \epsilon+ C_1 \epsilon^3 + C_2(M) \epsilon^{\kappa}  \|R \|_{H^{m}} ) \|R_2\|_{H^{m+1/2}}^2
\\ &\quad
+ 2 C_1 \epsilon^3  \|R_1 \|_{H^{m}}^2 + C_2(M) \epsilon^{\kappa}  \|R \|_{H^{m}}^3
+C_{res}\epsilon^{3}\left(1+\E(t)\right)
\end{align*}
under the assumption $ \|R\|_{H^{m}} \leq M $.

Suppose now
\begin{equation} \label{lenn1}
- c_1+ C_1 \epsilon^2 + C_1 \epsilon+ C_1 \epsilon^3 + C_2(M) \epsilon^{\kappa} M < 0 ,
\end{equation}
and
\begin{equation} \label{lenn2}
-  \frac{\eta}{2} +  C_1  +  C_1 \epsilon + C_2(M) \epsilon^{\kappa-2}  \|R \|_{H^{m}} < 0 .
\end{equation}
Then,  we end up with
\begin{align*}
 \partial_t \E & \leq
 2 C_1 \epsilon^3  \|R_1 \|_{H^{m}}^2 + C_2(M) \epsilon^{\kappa}  \|R \|_{H^{m}}^3
+C_{res}\epsilon^{3}\left(1+\E\right)
\\ & \leq 2 C_1 \epsilon^3  \E + \epsilon^3 \E +C_{res}\epsilon^{3}\left(1+\E\right),
\end{align*}
as long as
\begin{equation} \label{lenn3}
 C_2(M) \epsilon^{\kappa-3}  M  \leq 1 .
  \end{equation}
Under this condition,~\eqref{lenn2} simplifies into
  \begin{equation} \label{lenn4}
-  \frac{\eta}{2} +  2C_1  +  1 < 0 .
\end{equation}
We choose now $ \eta > 0 $ so large that~\eqref{lenn4} is satisfied.
Via Gr\"onwall's inequality we then deduce
\begin{align*} 
  \E(t) \leq C_{res} \epsilon^3 t e^{(2 C_1 + 1+ C_{res}) \epsilon^3 t} \leq C_{res}  \tau_1 e^{(2 C_1 + 1+ C_{res}) \tau_1} =: M^2,
\end{align*}
for $0 < \epsilon^3 t \leq \tau_1$. We are done if
we choose $ \varepsilon_0 > 0 $ so small that~\eqref{lenn1} and~\eqref{lenn3} are satisfied for all $ \epsilon \in (0,\epsilon_0) $
for this definition of $ M $.

\subsection{From Theorem~\ref{mainresult7} to Theorem~\ref{mainresult2}}
\label{sec78}

By proving Theorem~\ref{mainresult7}, we have obtained a KdV approximation result for the
$ \cZ $-variable. It remains to the transfer this result to the $ V $-variable, i.e.,
to conclude Theorem~\ref{mainresult2} from Theorem~\ref{mainresult7}.
We have
$$
\cZ(t) = S_{\omega(t)} S_{diag} S_{\theta} V(t) , \qquad
\textrm{resp.} \qquad V(t) =   S_{\theta}^{-1} S_{diag}^{-1}  S_{\omega(t)}^{-1}  \cZ(t).
$$
Combining Lemma~\ref{applem-1}, Lemma~\ref{applem-2}, and Lemma~\ref{applem0}
shows that $ S_{\omega(t)} S_{diag} S_{\theta} $ is an isomorphism
between $ H^\infty_{\mu(t),m-1}  \times H^\infty_{\mu(t),m} $ and  $ (H^m)^2 $ as long as $\mu(t) \geq 0$.
Since $ H^\infty_{\mu(t),m-1}  \times H^\infty_{\mu(t),m} $ can be embedded
in every $ H^s $-space for $0 \leq \epsilon^3 t \leq \tau_1$ the estimates~\eqref{impresolv} and~\eqref{Hmbound} follow immediately from
the estimates~\eqref{impresolvz}
 and~\eqref{Hmboundz}.

\section{Discussion}
\label{secdiscussion}

In~\S\ref{secbackU} we discuss the approximation result in the original coordinate system,
in~\S\ref{sobolev} what happens if  we work in Sobolev spaces
instead of spaces of analytic functions, in~\S\ref{secothers} whether other formally derived amplitude equations at the Eckhaus boundary are valid or not, and finally in~\S\ref{seccounterexample} we explain that although the KdV equation can also be derived  in the parameter region $ \A_h $, it makes wrong predictions in $\A_h$.

\subsection{The approximation result in the original variables}

\label{secbackU}

 The modulations $(\phi,s)$ in~\eqref{modsol} satisfy~\eqref{modeq2} in the coordinates~\eqref{newcoord}. To regain the phase $\phi$ from the local wave number  $\psi = \partial_x \phi$, we have to integrate the pointwise bound in Theorem~\ref{mainresult2}, leading to an error bound on a spatial interval of length $\mathcal{O}(\epsilon^{-\rho})$ with arbitrary -- but fixed -- $\rho > 0$, rather than an $\R$-uniform error bound.
 Moreover, we have to allow for a global phase $ e^{i \phi(X_0,T)} $.
These two restrictions  have been observed in other papers before, see
\cite[Section 4]{MS04b}
and
\cite[Corollary 2.2]{DS09JNS}. We consider the Ginzburg-Landau equation~\eqref{eq1}
in the original coordinates $ X,T $ and with some abuse of notation
we denote the variables depending on the new coordinates $ x,t$ introduced in~\eqref{newcoord}
by the same symbols.

We have to compare the modulated wave-train solution
$$
\Psi(X,T) =   \Psi_{\per}(X,T) \exp \left( s(X,T) + i  \int^{X}_{X_0}
\psi(X',T)d X'  + i \phi(X_0,T)\right)
$$
to~\eqref{eq1} with the approximation
$$
\Psi_{app}^{\kappa,\epsilon}(X,T) =  \Psi_{\per}(X,T) \exp \left(s_{app}^{\kappa,\epsilon}(X,T) + i  \int^{X}_{X_0}
    \psi_{app}^{\kappa,\epsilon}(X',T)d X'   \right),
$$
for some $X_0 \in \R$. By Theorem~\ref{mainresult2} it holds
\begin{align*}
\lefteqn{\left|\exp(- i \phi(X_0,T)) \Psi(X,T)-\Psi_{app}^{\kappa,\epsilon}(X,T)\right|} \\
& \leq C \left( \left|  s(X,T)-s_{app}^{\kappa,\epsilon}(X,T)
         \right|
       + \int^{X}_{X_0} \left| \psi( X', T)
       - \psi_{app}^{\kappa,\epsilon}(X',T)\right|  d X'\right) \\
& \leq C \left(\epsilon^{\kappa} + \int^{X}_{X_0} \epsilon^{\kappa}  d X'\right) \\
& \leq C \epsilon^{\kappa} (1+|X - X_0|) .
\end{align*}
Hence, using the improved approximations $
\Psi_{app}^{\kappa,\epsilon}
$ the approximation result holds uniformly on intervals
larger than the natural $ \mathcal{O}(1/\epsilon)$ spatial scale  of the KdV equation.
However,
due to
$\partial_T\phi(X_0,T)=\mathcal{O}(\epsilon^2)$, which implies
$\sup_{T\in[0,{\tau}_1/\epsilon^3]}|\phi(X_0,T)|=\mathcal{O}(1/\epsilon)$,
only the amplitude of the modulated wave train is well-approximated, but not
its position.
In spaces of sufficiently spatially localized functions the above estimates may be improved.

\subsection{From analytic initial conditions to Sobolev initial conditions?}
\label{sobolev}

It is a natural question whether we can replace the spaces $ H^{\infty}_{\mu,m} $
by classical  Sobolev spaces $ H^m $,
i.e., what happens if we give up analyticity in  a strip in the
complex plane and choose the solutions to the KdV equation to be only finitely many
times differentiable. We have no answer at this point and have to postpone the
question to future research.
The difficulty comes from the terms in the error equation associated to the nonlinear term
$  \gamma_{non} \partial_{\xi}(A^2) $ in the KdV equation, namely
$$
\varepsilon^2 {\theta} 2( B_1(\cW_1,R_1) + B_2(\cW_1,R_2))
$$
in~\eqref{expansion}, which is of  $ \mathcal{O}(\varepsilon^2) $, but its influence
has to be estimated $ \mathcal{O}(\varepsilon^3) $.
 Due to the marginal sideband instability, the smoothing of the linear part is too weak to gain additional powers of $ \epsilon $. This was the reason why we used the artificial smoothing through $ S_{\omega}(t) $.
If we do not want use the artificial smoothing, we have to find an energy
in which this term is $ \mathcal{O}(\varepsilon^3) $.

\subsection{Approximation results for other amplitude equations appearing
at the Eckhaus boundary}

\label{secothers}

As explained in the introduction, in~\cite{vH94} various other
amplitude equations for the description of slow modulations in time and
space of wave-train solutions to the Ginzburg-Landau equation have been derived near the Eckhaus boundary.

For $ \alpha = \beta $ the coefficients $ \gamma_{lin} $ and $ \gamma_{non} $ vanish
simultaneously and so the ansatz has to be modified to
\begin{align*}
\begin{split}
  \psi_{*,\epsilon}(x,t) &= \varepsilon^2 A(\varepsilon x , \varepsilon^4 t ) , \qquad s_{*,\epsilon}(x,t) = \varepsilon^2 B(\varepsilon x , \varepsilon^4 t ),
\end{split}
\end{align*}
and
a Cahn-Hilliard equation
\begin{equation*}
 \partial_{\tau} A  = \nu_2 \partial_{\xi}^2 A  + \nu_4 \partial_{\xi}^2 A + \nu_{non} \partial_{\xi}(A^2),
\end{equation*}
with real-valued coefficients $ \nu_2 $, $ \nu_4 < 0 $, and $ \nu_{non}  $
can be derived.  The justification of the Cahn-Hilliard approximation is again
a non-trivial task since solutions of $ \mathcal{O}(\epsilon^2) $
have to be estimated on an $ \mathcal{O}(1/\epsilon^4) $-time scale.
In case $ \alpha = \beta = 0 $ an approximation result has been established
in~\cite{Duell2}.
For $ \alpha = \beta \neq 0 $
we conjecture that the justification of the Cahn-Hilliard equation,
goes very similar to~\cite{Schn99,Duell}.

For $ (\alpha,\beta) $
in the parameter region $ \mathcal{A}_h $, the wave train~\eqref{equi1} becomes
unstable via a Hopf-Turing instability, cf. Figure~\ref{vanHarten}.
As  explained in~\cite{vH94} depending on the scaling a number of different
amplitude systems, such as a single Ginzburg-Landau equation,
 can be derived. The instability scenario is similar
to the one of pattern forming systems, with a conservation law, close to the
first instability of the homogeneous rest state. The validity of the Ginzburg-Landau approximation for such
pattern forming systems has been
considered in~\cite{HSZ11}. We expect
that the justification proof given in~\cite{HSZ11} can be adapted to cover the present
situation.

\subsection{Failure of the KdV approximation}
\label{seccounterexample}

\begin{figure}[h]
 \centering
  \begin{tikzpicture}[rotate=0,scale=1]
    \begin{axis}[
	xmin=-1.5, xmax=1.5,
	ymin=-0.3, ymax=0.5,
	axis lines=center,
	ticks=none,
	ylabel={$\textrm{Re}(\lambda_+)$},
	xlabel={$k$},
	width=0.6\textwidth,
    ]
      \addplot+[black,thick, domain=-0.3:0.3,samples=100,no marks,opacity=0.5] {-3*x^2};
      \addplot+[black,thick, domain=0.6:1.4,samples=100,no marks,opacity=0.5] {-3*(x-1)^2+0.3};
      \addplot+[black,thick, domain=-1.4:-0.6,samples=100,no marks,opacity=0.5] {-3*(x+1)^2+0.3};
      \addplot[no marks] coordinates {(-0.05,0.3) (0.05,0.3)};
      \addplot[no marks] coordinates {(1,-0.025) (1,0.025)};
      \node[below,black] at (axis cs:0.96,-0.03) {$ \mathcal{O}(1) $};
      \node[left=4pt] at (axis cs:0.05,0.3) {$  \mathcal{O}(1) $};
    \end{axis}
  \end{tikzpicture}
  \begin{tikzpicture}[rotate=0,scale=1]
    \begin{axis}[
	xmin=-1.5, xmax=1.5,
	ymin=-0.3, ymax=0.5,
	axis lines=center,
	ticks=none,
	ylabel={$\textrm{Re}(\lambda_+)$},
	xlabel={$k$},
	width=0.6\textwidth,
    ]
      \addplot+[black,thick, domain=-0.5:0.5,samples=100,no marks,opacity=0.5] {-5*x^2*(x^2-0.2)};
      \addplot+[black,thick, domain=0.6:1.4,samples=100,no marks,opacity=0.5] {-3*(x-1)^2+0.3};
      \addplot+[black,thick, domain=-1.4:-0.6,samples=100,no marks,opacity=0.5] {-3*(x+1)^2+0.3};
      \addplot[no marks] coordinates {(-0.05,0.3) (0.05,0.3)};
      \addplot[no marks] coordinates {(-0.05,0.05) (0.05,0.05)};
      \addplot[no marks] coordinates {(1,-0.025) (1,0.025)};
      \node[below,black] at (axis cs:0.96,-0.03) {$ \mathcal{O}(1) $};
      \node[left=4pt] at (axis cs:0.05,0.08) {$ \mathcal{O}(\varepsilon^4) $};
      \node[left=4pt] at (axis cs:0.05,0.3) {$  \mathcal{O}(1) $};
    \end{axis}
  \end{tikzpicture}
 \caption{The spectral curve $ \lambda_+ $   in the parameter region $ \mathcal{A}_h $
before and   after the sideband instability occurred for a wave train which is already
Hopf-Turing unstable. }
 \label{Failurefig}
\end{figure}
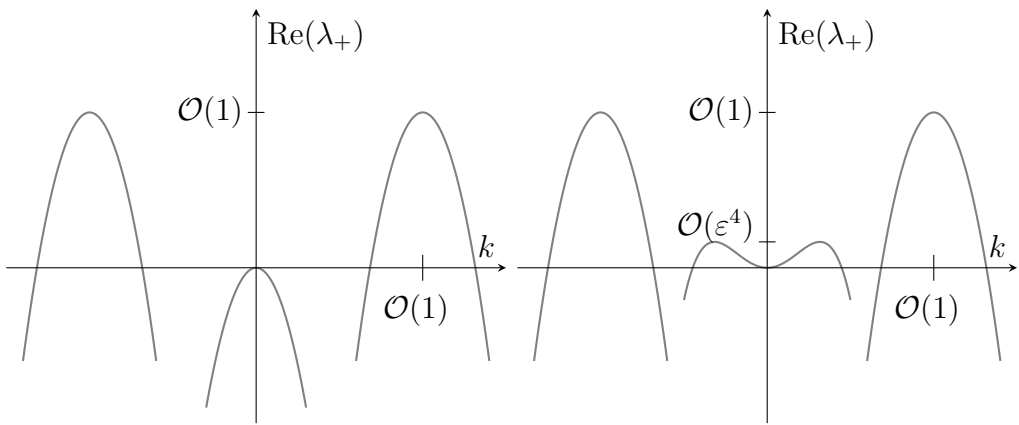

Due to its long-wave character,
also in the region $ \mathcal{A}_h $, see Figure~\ref{vanHarten},  a KdV equation can be formally derived,
although there is an
$ \mathcal{O}(1) $-instability
at a wave number $ k_{1} \neq 0$, cf. Figure~\ref{Failurefig}.
Thus, the derivation of the KdV equation in~\S\ref{sec3} is independent of whether we are in
 the parameter region $ \mathcal{A}_s $ or   $ \mathcal{A}_h $. However,
 the justification analysis is not.

 For solutions to the KdV equation  in $ H^\infty_{\mu_{\A},m} $
 the corresponding solution to the Ginzburg-Landau equation~\eqref{eq1}
is initially of order $ \mathcal{O}(e^{-\mu |k|/\varepsilon}) $ at a wave number $ k $.
For a wave number $ k_1 $ with
$ \mathcal{O}(1) $ growth rates at $ k_1 $, cf. Figure~\ref{Failurefig}, already for $ t = \mathcal{O}(1/\varepsilon) $
the solution is $ \mathcal{O}(1) $ and for $ t = \mathcal{O}(1/\varepsilon^3) $ the growth rate
is $ \mathcal{O}(e^{1/\varepsilon^2}) $. Therefore, the original system is expected to behave completely different as
predicted by the KdV equation.  For $ 2\pi/k_1 $-spatially periodic
solutions this  can possibly be made  rigorous
with the help of a center manifold reduction similar to~\cite{Schn96MN}.

\appendix

\section{Some technical estimates}

\subsection{Estimates for the change of coordinates}
\label{secsdiag}

We compute
\begin{align*} \hat{S_{diag}}(k) &= \begin{pmatrix} 1& -\frac{(1 - \upsilon(k))\hat{\theta}(k)}{\gamma(k)}\\ -1 & \frac{(1 + \upsilon(k))\hat{\theta}(k)}{\gamma(k)}\end{pmatrix} =:  \begin{pmatrix}
  s_{11}(k) &  \hat{\theta}(k) s_{12}(k)\\
  s_{21}(k) & s_{22}(k)
 \end{pmatrix}, \\ \qquad \hat{S_{diag}}(k)^{-1} &=
\frac{1}{2 \upsilon(k)}
\begin{pmatrix}
1 + \upsilon(k) &  1 - \upsilon(k)\\
\hat{\theta}(k)^{-1} \gamma(k) & \hat{\theta}(k)^{-1}\gamma(k)
  \end{pmatrix},\end{align*}
where we recall from~\S\ref{sec2} that $\gamma(k) = (\alpha k^2 - 2 i k)/\sigma$ and $ \upsilon(k)$ is the principal square root $\sqrt{1-\gamma(k)^2-2 \beta \gamma(k)}$. One readily verifies $s_{ij} \in L^\infty(\R,\C)$ for $i,j \in \{1,2\}$ and $\hat{S_{diag}},\hat{S_{diag}}^{-1} \in L^\infty(\R,\mathrm{Mat}_{2 \times 2}(\C))$.
As a direct consequence of these pointwise estimates we find
\begin{lemma}\label{applem-1}
For
 $\mu \ge 0$ and $s \geq 0$
the linear mapping  $ S_{diag} :H_{\mu,s}^{\infty} \to H_{\mu,s}^{\infty} $
 is  bijective and bounded  with bounded inverse.
 The similar statement is true in the $ W_{\mu,m} $-spaces.
\end{lemma}
\noindent
{\bf Proof.}
We have
$$
\| S_{diag}^{\pm 1}  u \|_{H_{\mu,s}^{\infty}} \leq \| \hat{S_{diag}}^{\pm 1} \|_{L^\infty(\R,\mathrm{Mat}_{2 \times 2}(\C))}  \| u \|_{H_{\mu,s}^{\infty}} \leq C_{diag}
 \| u \|_{H_{\mu,s}^{\infty}} .
$$
The estimates in $ W_{\mu,s} $ follow similarly.
\qed

Analogously, we establish
\begin{lemma}\label{applem-2}
For
 $\mu \ge 0$ and $s \geq 0$
the linear mapping $ S_{\theta} :H_{\mu,s-1}^{\infty} \times H_{\mu,s}^{\infty}  \to (H_{\mu,s}^{\infty})^2 $ is  bijective and bounded  with bounded inverse.
\end{lemma}

Finally we have
\begin{lemma}\label{applem0}
For $t\in[0,\mu_*/( \eta \epsilon^{3})]$ the linear mapping
$ S_{\omega}(t)  : H^\infty_{\mu(t),s}  \to  H^{s} $,
with  $ \mu(t) = (\mu_* - \eta \epsilon^3 t)/\epsilon $ is bijective and bounded  with bounded inverse.
The similar statement is true in the $ W_{\mu,s} $-spaces.
\end{lemma}
\noindent
{\bf Proof.}
This follows immediately from the definitions.
\qed

\subsection{Estimates for the nonlinear terms}
\label{appA2}

\begin{lemma}\label{applem1}
The spaces $H_{\mu,s}^{\infty}$ are Banach algebras for
 $\mu \ge 0$ and $s > \frac{1}{2}$.
 In detail, there exists a $ \mu $-independent constant $ C_{s} $ such that
 $$
 \| u v \|_{H_{\mu,s}^{\infty}} \leq  C_{s}  \| u \|_{H_{\mu,s}^{\infty}}  \| v \|_{H_{\mu,s}^{\infty}}
 $$
 for all $ u,v \in H_{\mu,s}^{\infty}$.
\end{lemma}
\noindent
\begin{proof}
 Suppose that $u$ and $v$ are in $H_{\mu,s}^{\infty}$.  Then
 \begin{align*}
  \lVert u v \rVert_{H_{\mu,s}^{\infty}} &\leq \lVert (1+k^2)^{s/2} (|\e^{\mu |\cdot|}\hat u|\ast|\e^{\mu |\cdot|}\hat v| )(k) \rVert_{L^2(k)} \\
  &\leq \tilde C_s \left(\lVert (1+k^2)^{s/2} \e^{\mu |k|}\hat u(k) \rVert_{L^2(k)} \lVert \e^{\mu |k|}\hat v(k) \rVert_{L^1(k)} \right.
  \\&\qquad \left.+ \lVert \e^{\mu |k|}\hat u(k) \rVert_{L^1(k)} \lVert (1+k^2)^{s/2} \e^{\mu |k|}\hat v(k) \rVert_{L^2(k)}\right)\\
  &\leq C_s \lVert (1+k^2)^{s/2} \e^{\mu |k|}\hat u(k) \rVert_{L^2(k)} \lVert (1+k^2)^{s/2}\e^{\mu |k|}\hat v(k) \rVert_{L^2(k)} \\
  &= C_s\lVert u \rVert_{H_{\mu,s}^{\infty}} \lVert v \rVert_{H_{\mu,s}^{\infty}},
 \end{align*}
 where we used the continuous embedding ${L^2(\R,(1+|x|^2)\mathrm{d}x)} \hookrightarrow L^1(\R)$ for $s > \frac{1}{2}$ and Young's inequality.
\end{proof}

%

For  our error estimates we need the following tame estimates.
\begin{corollary}\label{appcoro1}
For $ \delta > 0 $,
 $\mu \ge 0$ and $s > 1/2 $ we have
 $$
 \| u^2 \|_{H_{\mu,s}^{\infty}} \leq  C_{s}  \| u \|_{H_{\mu,1/2+\delta}^{\infty}}  \| u \|_{H_{\mu,s}^{\infty}}
 $$
 for all $ u \in H_{\mu,s}^{\infty}$.
\end{corollary}
\noindent
\begin{proof}
 This follows obviously from the proof of Lemma~\ref{applem1}.
\end{proof}


\begin{corollary}
\label{applem5}
For
 $\mu \ge 0$ and $s \geq 0 $ we have
 $$
 \| u v  \|_{H_{\mu,s}^{\infty}} \leq  C_{s}  \| u \|_{W_{\mu,s}}  \| v \|_{H_{\mu,s}^{\infty}}
 $$
 for all $ u  \in W_{\mu,s} $ and $ v\in H_{\mu,s}^{\infty}$.
\end{corollary}
\noindent
\begin{proof}
 This follows obviously from the proof of Lemma~\ref{applem1}.
\end{proof}


\begin{corollary}\label{applem2}
The entire function $ h: \C \to \C $ defined by $h(s) = e^{2s} - 1 - 2s = \sum_{j = 2}^\infty h_js^j$ with $ h_j = 2^j/j! $ defines via
$ u \mapsto h(u) $ an entire function in $H_{\mu,s}^{\infty}$  for
 $\mu \ge 0$ and $s > \frac{1}{2}$.
\end{corollary}
\noindent
\begin{proof}
 The series converges absolutely, since $h$ is an entire function, $h(0) = 0$, and $H_{\mu,s}^{\infty}$ is a Banach algebra for $\mu \geq 0$ and $s > \frac{1}{2}$ according to Lemma~\ref{applem1}.
\end{proof}
%
%

\begin{lemma}
\label{applem3}
The nonlinear mapping $N_{\cZ}(t)$ is smooth
from $ H^{s+1} $ to $ H^{s} $ for $ s \geq 2 $ and $t \in [0,\mu_*/(\eta\epsilon^3)]$.
\end{lemma}
\noindent
{\bf Proof.}
From the explicit representation~\eqref{nonlN1} of $ N_Y $,  Lemma~\ref{applem1}
and Corollary~\ref{applem2}  it follows that $ N_Y $ is a smooth mapping
from $H_{\mu,s+1}^{\infty}$ to
$H_{\mu,s}^{\infty}$  for $\mu \ge 0$ and $s > \frac{3}{2}$. From Lemma
\ref{applem0} it follows that $ N_\Y $ is a smooth mapping $ H^{s+1} $ to $ H^{s} $.
Finally, Lemma
\ref{applem-1} implies that $ N_{\cZ} $ is a smooth mapping $ H^{s+1} $ to $ H^{s} $ (see also Figure~\ref{figvariables}).
\qed

\subsection{An inequality  for the residual estimates}

The formal expansion of the curve of eigenvalues and of the kernels
in the multilinear maps can be estimated with the aid of the following lemma.
\begin{lemma} \label{applem4}
Let $  {\theta}_0\geq 0 $, $ {\theta}_{\infty} \in \R $, and let $ g:\R \rightarrow \C $ satisfy
$$ | g(k) | \leq C \min (| k |^{{\theta}_0}, (1+|k|)^{{\theta}_{\infty}}).$$
Then for the associated multiplication operator $ g_{op} =  \F^{-1} g \F $
the following holds. For i) $ \mu_1 > \mu_2 $ and $ m_1,m_2 \geq 0 $ or ii) $ \mu_1 = \mu_2 $ and $ m_2 -m_1 \geq
\max({\theta}_0, {\theta}_{\infty}) $
we have
$$
\| g_{op} A(\epsilon \cdot ) \|_{H^{\infty}_{\mu_1/\epsilon,m_1}}
\leq C \varepsilon^{{\theta}_0-1/2} \| A( \cdot ) \|_{H^{\infty}_{\mu_2,m_2}}
$$
for all $ \varepsilon \in (0,1) $.
\end{lemma}
\noindent
{\bf Proof.}
This follows immediately from the fact that the
 left-hand side of this inequality can be estimated by
\begin{align*}
& \leq   \sup_{k \in \R} \left| g(k)  \left(1 + \left(\frac{k}{\epsilon}\right)^2\right)^{\frac{m_1-m_2}{2}}
e^{\left(\mu_1-\mu_2\right)|k|/\epsilon}
\right|
\
 \|  A (\epsilon \cdot )  \|_{H^{\infty}_{\mu_2/\epsilon,m_2}}
 \\
&  \leq
C \epsilon^{{\theta}_0}
\  \epsilon^{-1/2}  \|  A ( \cdot )  \|_{H^{\infty}_{\mu_2,m_2}}
\end{align*}
where the loss of $ \epsilon^{-1/2} $ is  due to  the scaling properties of the $ L^2 $-norm.
\qed

In $ W_{\mu,m} $-spaces there is no $ \epsilon^{-1/2} $ loss due to the scaling invariance of the norm and so we have as a direct consequence:
\begin{corollary} \label{appcoro4}
Let $  {\theta}_0\geq 0 $, $ {\theta}_{\infty} \in \R $, and let
$ g(k) $ satisfy
$$ | g(k) | \leq C \min (| k |^{{\theta}_0}, (1+|k|)^{{\theta}_{\infty}}).$$
Then for the associated operator $ g_{op} =  \F^{-1} g \F $
the following holds. For i) $ \mu_1 > \mu_2 $ and $ m_1,m_2 \geq 0 $ or ii) $ \mu_1 = \mu_2 $ and $ m_2 -m_1 \geq
\max({\theta}_0, {\theta}_{\infty}) $
we have
$$
\| g_{op} A(\epsilon \cdot ) \|_{W_{\mu_1/\epsilon,m_1}}
\leq C \varepsilon^{{\theta}_0} \| A( \cdot ) \|_{W_{\mu_2,m_2}}
$$
for all $ \varepsilon \in (0,1) $.
\end{corollary}

\bibliographystyle{plain} 
\bibliography{GLbib12}
\small

\end{document}